\documentclass[12pt, a4paper]{amsart}

\usepackage{geometry}
\usepackage{amsmath, amssymb}
\usepackage{enumerate}
\usepackage{diagbox}
\usepackage{graphicx}
\usepackage{color}
\graphicspath{{fig/}}

\newtheorem{example}{Example}[section]
\newtheorem{remark}{Remark}
\newtheorem{theorem}{Theorem}
\newtheorem{algorithm}{Algorithm}

\makeatletter %
\@addtoreset{equation}{section} %
\makeatother

\DeclareMathOperator*{\argmin}{argmin}

\def\half{\frac 1 2}

\newcommand{\vr}{\vec{r}}

\newcommand{\mbd}{\mathbb{D}}
\newcommand{\mbfd}{\mathbf{d}}
\newcommand{\mbe}{\mathbb{E}}

\newcommand{\mcf}{\mathcal{F}}
\newcommand{\mcg}{\mathcal{G}}
\newcommand{\mcm}{\mathcal{M}}
\newcommand{\mcs}{\mathcal{S}}
\newcommand{\mcw}{\mathcal{W}}

\newcommand{\br}{\mathbb{R}}
\newcommand{\rd}{\mathrm{d}}
\newcommand{\veps}{\varepsilon}

\newcommand{\bd}{\boldsymbol{d}}
\newcommand{\bn}{\boldsymbol{n}}
\newcommand{\boeta}{\boldsymbol{\eta}}
\newcommand{\bx}{\boldsymbol{x}}
\newcommand{\bxi}{\boldsymbol{\xi}}

\newcommand{\tf}{\tilde{f}}
\newcommand{\tg}{\tilde{g}}
\newcommand{\tr}{\tilde{r}}

\newcommand{\abs}[1]{\left\vert#1\right\vert}
\newcommand{\brac}[1]{\left(#1\right)}
\newcommand{\innprod}[1]{\langle#1\rangle}
\newcommand{\norm}[1]{\left\Vert#1\right\Vert}

\newcommand{\nn}{\nonumber}

\begin{document}

\title[The quadratic Wasserstein metric for earthquake location]{The quadratic Wasserstein metric for earthquake location}
  
\author{Jing Chen}
\address{Department of Mathematical Sciences \\
  Tsinghua University \\
  Beijing, China 100084 \\
  email: jing-che16@mails.tsinghua.edu.cn }
  
\author{Yifan Chen}
\address{Department of Mathematical Sciences \\
  Tsinghua University \\
  Beijing, China 100084 \\
  email: chenyifan14@mails.tsinghua.edu.cn }

\author{Hao Wu$^*$}
\address{$^*$Corresponding author \\
  Department of Mathematical Sciences \\
  Tsinghua University \\
  Beijing, China 100084 \\
  email: hwu@tsinghua.edu.cn }

\author{Dinghui Yang}
\address{Department of Mathematical Sciences \\
  Tsinghua University \\
  Beijing, China 100084 \\
  email: dhyang@math.tsinghua.edu.cn }

\date{\today}

\begin{abstract}
	In [Engquist et al., Commun. Math. Sci., 14(2016)], the Wasserstein metric was successfully introduced to the full waveform inversion. We apply this method to the earthquake location problem. For this problem, the seismic stations are far from each other. Thus, the trace by trace comparison [Yang et al., arXiv(2016)] is a natural way to compare the earthquake signals.  
	
	Under this framework, we have derived a concise analytic expression of the Fr\`echet gradient of the Wasserstein metric, which leads to a simple and efficient implementation for the adjoint method. We square and normalize the earthquake signals for comparison so that the convexity of the misfit function with respect to earthquake hypocenter and origin time can be observed numerically. To reduce the impact of noise, which can not offset each other after squaring the signals, a new control parameter is introduced. Finally, the LMF (Levenberg-Marquardt-Fletcher) method is applied to solve the resulted optimization problem. According to the numerical experiments, only a few iterations are required to converge to the real earthquake hypocenter and origin time. Even for data with noise, we can obtain reasonable and convergent numerical results.
	 \medskip
	
	\noindent {\bf Keywords:} Optimal transport, Wasserstein metic, Inverse theory, Waveform inversion, Earthquake location
\end{abstract}

\maketitle


\section{Introduction} \label{sec:intro}

The Wasserstein metric is an important concept in the optimal transport theory \cite{Sa:15, Vi:03, Vi:08}. It measures the difference between two probability distributions as the optimal cost of rearranging one distribution into the other. This kind of problem was first proposed by French engineer Gaspard Monge \cite{Mo:81}. He wanted to find a way to move a pile of sand to a designated location at a minimum cost. As the metric provides a global comparison tool, it is very suitable to model and solve problems from computer vision \cite{RuToGu:98}, machine learning \cite{ArChBo:17}, etc.

From the mathematical point of view, there are many advantages of the Wasserstein metric \cite{AmGi:13, EnFr:14, EnFrYa:16, Vi:03}, especially for the quadratic Wasserstein metric ($W_2$), e.g. the convexity with respect to shift, dilation and partial amplitude change, the insensitivity with respect to noise. In \cite{EnFr:14}, Engquist and Froese first used this metric to measure the misfit between seismic signals. The idea was then developed to invert the velocity structure \cite{EnFrYa:16}. Due to the convexity property of the Wasserstein metric, the full waveform inversion converges to the correct solution, even from poor initial values. In \cite{YaEnSuFr:16}, the method was further applied to more realistic examples. Motivated by this idea, M\'etivier and collaborators proposed the KR norm based full waveform inversion \cite{MeBrMeOuVi:16, MeBrMeOuVi:16b}. They also show the superiority of their method through some realistic examples. Different from Engquist's method, the KR norm is related to the Wasserstein metric with linear cost function.

In this study, we would like to apply the quadratic Wasserstein metric to the earthquake location problem. The earthquake location is a fundamental problem in seismology \cite{Ge:03, Ma:15, Th:14}. In quantitative seismology, there are many applications \cite{LeSt:81, PrGe:88, SaLoZo:08, Th:85}. Traditionally, the ray theory based earthquake location methods have been widely used \cite{Ge:03, Ge:03b, Ge:12, WaEl:00}. But this method is of low accuracy when the seismic wave length is not small enough compared to the scale of wave propagation region \cite{EnRu:03, JiWuYa:08, RaPoFi:10, WuYa:13}. Therefore, the waveform based earthquake location methods \cite{KiLiTr:11, LiPoKoTr:04, ToYaLiYaHa:16, WuChHuYa:16, WuChJiJiYa:17} are developed to determine the earthquakes' parameters accurately. However, it is well known that the waveform inversion with $\ell^2$ norm is suffering from the cycle-skipping problem. It requires accurate initial data for inversion. For the earthquake location problem, the situation is even worse, since the seismic focus is modeled by the highly singular delta function $\delta(\bx-\bxi)$ \cite{AkRi:80, Ma:15}. Fortunately, the Wasserstein metric has been shown to be effective in overcoming the famous cycle-skipping problem \cite{EnFrYa:16, MeBrMeOuVi:16, MeBrMeOuVi:16b, YaEnSuFr:16}. This is the motivation for us to work on this topic.

For the earthquake location problems, the receivers are located far apart from each other. This is different from the situation that the receivers are close to each other in the problems of exploration seismology. Thus, we would like to follow the idea of trace by trace comparison with $W_2$ metric \cite{YaEnSuFr:16}. In this paper, we apply this metric to invert the seismic focus parameters. A concise analytic expression of the Fr\`echet gradient is derived to  simplify the numerical computation. We use the LMF method \cite{Fl:91, Le:44, MaNiTi:04, Ma:63} to solve the optimization problem since it has a least square structure. The numerical experiments show that the computational efficiency is greatly improved. We also want to point out that the method developed in this paper may be effective for geological scale problems.

This paper is organized as follows. After reviewing the formulation and basic properties of the $W_2$ metric in Section \ref{sec:qwm}, we apply it to the earthquake location problem in Section \ref{sec:ael}. Since we are considering the trace by trace comparison of the $W_2$ metric, it is very easy to derive the Fr\'echet gradient and the sensitivity kernel. The efficient LMF method is introduced to solve the optimization problem. In Section \ref{sec:num}, the numerical experiments are provided to demonstrate the effectiveness and efficiency of the method. Finally, we conclude the paper in Section \ref{sec:con}.

\section{The quadratic Wasserstein metric} \label{sec:qwm}
Let $\tf$ and $\tg$ be two probability density functions on $\br$, then the mathematical definition of the quadratic Wasserstein metric between $\tf$ and $\tg$ is formulated as follows \cite{EnFr:14, EnFrYa:16}:
\begin{equation} \label{eqn:qwm}
	W_2^2(\tf,\tg)=\inf_{T\in\mcm}\int_{\br}\abs{t-T(t)}^2\tf(t)\rd t,
\end{equation}
in which $\mcm$ is the set of all the rearrange maps from $\tf$ to $\tg$. According to the ``Optimal transportation theorem for a quadratic cost on $\br$'' (see P74 in \cite{Vi:03}), the optimal transportation cost and the optimal map are
\begin{equation} \label{eqn:qwm_sol}
	W_2^2(\tf,\tg)=\int_0^1\abs{F^{-1}(t)-G^{-1}(t)}^2\rd t, \quad T(t)=G^{-1}(F(t)),
\end{equation}
where $F(t)$ and $G(t)$ are the corresponding cumulative distribution functions of $\tf(t)$ and $\tg(t)$:
\begin{equation*}
	F(t)=\int_{-\infty}^t \tf(\tau)\rd \tau, \quad G(t)=\int_{-\infty}^t \tg(\tau)\rd \tau.
\end{equation*}
It is easy to verify that
\begin{equation} \label{eqn:qwm_dom}
	\tf(t)=\tg(T(t))T'(t).
\end{equation}
In this study, we prefer to use the 1D quadratic Wasserstein metric. This is based on the fact that the receivers are far from each other in geological scale problems. Moreover, the solution of the 1D quadratic Wasserstein metric can be easily obtained by \eqref{eqn:qwm_sol}.

In seismology, there are difficulties in applying the Wasserstein metric. Firstly, the seismic signal $f(t)$ and $g(t)$ are not positive. For example, in many situations the seismogram at source has the form of Ricker wavelet
\begin{equation*}
	R(t)=A(1-2\pi^2f_0^2t^2)e^{-\pi^2f_0^2t^2}.
\end{equation*}
Here $f_0$ is the dominant frequency and $A$ is the normalization factor. This Ricket wavelet is not always positive over the entire time axis. Secondly, the comparison between $\tf$ and $\tg$ under the Wasserstein metric requires the mass conservation, i.e.
\begin{equation*}
	\int_{\br}\tf(t)\rd t=\int_{\br}\tg(t)\rd t.
\end{equation*}
The above mentioned difficulties can be easily solved by considering the following reformulated distance
\begin{equation} \label{eqn:ref_dis}
	\mbfd(f,g)=W_2^2\brac{\frac{f^2}{\innprod{f^2}},\frac{g^2}{\innprod{g^2}}},
\end{equation}
in which, the operator $\innprod{\cdot}$ denote the integral over the real axis
\begin{equation*}
	\innprod{f}=\int_{\br}f(t)\rd t.
\end{equation*}

\begin{remark}
	In Section 3.1 of manuscript \cite{YaEnSuFr:16}, the authors prefer to add a constant $c$ to ensure the positivity. However, in our numerical tests, the square strategy seems to be more suitable for the earthquake location problems.
\end{remark}

\begin{remark}
	In \cite{EnFrYa:16}, the convexity of the quadratic Wasserstein metric with respect to shift, stretching and partial amplitude loss has been proved. Thus, we will not repeat here. 
\end{remark}

The aforementioned discussion is concerning the theoretical model. For practical problems in seismology, the signals $f(t)$ and $g(t)$ can be considered to have compact support $[0,t_f]$ for $t_f$ large enough. Thus, the operator $\innprod{\cdot}$ is redefined as
\begin{equation*}
	\innprod{f(t)}=\int_0^{t_f}f(t)\rd t.
\end{equation*}
In the later part of the paper, we will default to this notation unless otherwise specified.

\subsection{The Fr\'echet gradient}
We have defined the distance $\mbfd(f,g)$ based on Wasserstein metric. To solve the resulted optimization problem, it is necessary to derive the Fr\'echet gradient $\nabla_f \mbfd$. Define
\begin{equation*}
	\mcs(f)=\frac{f^2}{\innprod{f^2}}, \quad \mcw(\mcf,\mcg)=W_2^2(\mcf,\mcg),
\end{equation*}
then we have
\begin{equation*}
	\mbfd(f,g)=\mcw(\mcs(f),\mcs(g)),
\end{equation*}
and
\begin{equation} \label{eqn:chain_rule}
	\nabla_f \mbfd=\nabla_{\mcf} \mcw \cdot\nabla_f\mcs.
\end{equation}

Before derivation, we have to emphasize that all the high order terms are ignored without any explanation. We first derive the gradient $\nabla_{\mcf}\mcw$. Let $\delta \mcf$ be a small perturbation\footnote{In order to avoid repeating the explanation, we use $\delta f$ to denote the small perturbation of arbitrary function $f$ in the later part of the paper.} of $\mcf$, according to \eqref{eqn:qwm}-\eqref{eqn:qwm_sol}
\begin{equation*}
	\mcw+\delta \mcw=\int_0^{t_f}\abs{t-(T+\delta T)}^2(\mcf+\delta \mcf)\rd t.
\end{equation*}
This leads to
\begin{equation} \label{eqn:delta_mcw}
	\delta \mcw=\int_0^{t_f}\abs{t-T}^2\delta \mcf\rd t-2\int_0^{t_f}(t-T)\mcf \delta T\rd t.
\end{equation}
Since the Wasserstein metric measures the difference between two probability density functions, we can naturally assume that
\begin{equation} \label{asm:deltamcf}
	\int_0^{t_f}\delta \mcf(t)\rd t=0,
\end{equation}
Using the equation \eqref{eqn:qwm_dom}, we get
\begin{equation*}
	\mcf+\delta \mcf=\brac{\mcg(T)+\mcg'(T)\delta T}(T'+(\delta T)'),
\end{equation*}
which yields
\begin{equation*}
	\delta \mcf=\mcg(T)(\delta T)'+\mcg'(T)T'\delta T=(\mcg(T)\delta T)'.
\end{equation*}
Integrating the above equation over $[0,t]\subset[0,t_f]$, leads to
\begin{equation} \label{eqn:delta_T}
	\mcg(T(t))\delta T(t)=\int_0^t \delta \mcf(\tau)\rd \tau+\mcg(T(0))\delta T(0)=\int_0^t \delta \mcf(\tau)\rd \tau,
\end{equation}
where the second equality holds since $T(0)=0$ and
\begin{equation*}
	\mcg(T(0))=\mcg(0)=\mcs(g(0))=0.
\end{equation*}
Using \eqref{eqn:delta_T} in \eqref{eqn:delta_mcw}, we have, in light of \eqref{eqn:qwm_dom}
\begin{multline} \label{eqn:dwdf}
	\delta \mcw=\int_0^{t_f}\abs{t-T(t)}^2\delta \mcf \rd t-2\int_0^{t_f} (t-T(t))T'(t) \int_0^t\delta \mcf(\tau)\rd \tau\rd t \\
		=\int_0^{t_f}\abs{t-T(t)}^2\delta \mcf \rd t-2\int_0^{t_f} \brac{\int_t^{t_f}(\tau-T(\tau))T'(\tau)\rd \tau} \delta \mcf(t)\rd t
		=\int_0^{t_f}\varphi(t)\delta \mcf(t)\rd t.
\end{multline}
where
\begin{equation*}
	\varphi(t)=\abs{t-T(t)}^2-2\int_t^{t_f}(\tau-T(\tau))T'(\tau)\rd \tau+C,
\end{equation*}
here $C\in\br$ is a constant. Differentiating both sides of the above equation with respect to $t$ gives
\begin{equation*}
	\varphi'(t)=2(t-T(t)).
\end{equation*}
Thus, the simplest form of the function $\varphi(t)$ can be written as
\begin{equation} \label{eqn:varphi}
	\varphi(t)=2\int_0^t(\tau-T(\tau))\rd \tau.
\end{equation}
\begin{remark}
	We can also get this expression by the optimal transport theory and the duality theory of linear programming. For interested readers, we refer to P14 Theorem 1.17 in \cite{Sa:15} for the details.
\end{remark}

Next, we would like to derive $\nabla_{f}\mcs$. For small perturbation $\delta f$ of $f$, we have
\begin{equation*}
	\mcs+\delta \mcs=\frac{(f+\delta f)^2}{\innprod{(f+\delta f)^2}}=\frac{f^2+2f\delta f}{\innprod{f^2+2f\delta f}}
		=\brac{f^2+2f\delta f}\brac{\frac{1}{\innprod{f^2}}-\frac{2\innprod{f\delta f}}{\innprod{f^2}\innprod{f^2+2f\delta f}}}.
\end{equation*}
It follows that
\begin{equation} \label{eqn:dsdf}
	\delta \mcs=\frac{2f\delta f}{\innprod{f^2}}-\frac{2f^2\innprod{f\delta f}}{\innprod{f^2}^2}.
\end{equation}

In summary, the Fr\'echet gradient can be obtained by combining \eqref{eqn:chain_rule} and \eqref{eqn:dwdf}-\eqref{eqn:dsdf}
\begin{multline} \label{eqn:frechet_gradient}
	\delta \mbfd=\int_0^{t_f}\brac{2\int_0^t(\tau-T(\tau))\rd \tau}\brac{\frac{2f\delta f}{\innprod{f^2}}-\frac{2f^2\innprod{f\delta f}}{\innprod{f^2}^2}}\rd t \\
		=\int_0^{t_f}4\brac{A(t)-B}f(t)\delta f(t)\rd t,
\end{multline}
where
\begin{equation*}
	A(t) = \frac{\int_0^t(\tau-T(\tau))\rd \tau}{\int_0^{t_f}f^2(t)\rd t}, \quad
	B = \frac{\int_0^{t_f}\brac{\int_0^t(\tau-T(\tau))\rd \tau}f^2(t)\rd t}{\brac{\int_0^{t_f}f^2(t)\rd t}^2}.
\end{equation*}

\subsection{Reduce the impact of noise}
We now turn to discuss the impact of noise on the reformulated distance $\mbfd(f,g)$. In \eqref{eqn:ref_dis}, we take the square of the signals $f,\;g$ to ensure the positivity, which implies that Theorem 3.1 in \cite{EnFrYa:16} is not directly applicable here. In the following, we consider a more general situation.

\begin{theorem} \label{thm:noise1}
	Let $\tg(t):[0,1]\to (0,M_1]$ and
	\begin{align*}
		& \tf_N(t)=\tg(t)+\tr_N(t), \\
		& \tr_N(t)=\left\{\begin{array}{ll}
			\tr_1, & t\in[0,\frac{1}{N}], \\
			\tr_2, & t\in(\frac{1}{N},\frac{2}{N}], \\
			\cdots \\
			\tr_N, & t\in(\frac{N-1}{N},1],
		\end{array}\right.
	\end{align*}
	in which $\tr_j$ are i.i.d. random variables with zero expectation and bounded variance
	\begin{equation*}
		\mbe \tr_j=0, \quad \mbd \tr_j<+\infty, \quad j=1,2,\cdots,N.
	\end{equation*}
	We further assume that $\tf_N(t):[0,1]\to (0,M_2]$, then $\mbe W_2^2(\tf_N/\innprod{\tf_N},\tg/\innprod{\tg})=O(1/N)$.
\end{theorem}
The proofs are almost identical to Theorem 3.1 in \cite{EnFrYa:16}. Thus, they will not be reproduced here.

For practical problems, consider the signal $f_N(t),\;g(t)$ defined on $[0,t_f]$ and
\begin{equation} \label{eqn:fNg}
	f_N(t)=g(t)+r_N(t).
\end{equation}
Here
\begin{equation} \label{eqn:rNt}
	r_N(t)=\left\{\begin{array}{ll}
		r_1, & t\in[0,\frac{t_f}{N}], \\
		r_2, & t\in(\frac{t_f}{N},\frac{2t_f}{N}], \\
		\cdots \\
		r_N, & t\in(\frac{(N-1)t_f}{N},t_f],
	\end{array}\right.
\end{equation}
in which $r_j$ are i.i.d. random variables with bounded expectation and variance
\begin{equation} \label{eqn:iidrv}
	\mbe r_j=\mu, \quad \mbd r_j=\sigma^2, \quad j=1,2,\cdots,N.
\end{equation}
Taking account of noise, we redefine the distance function in \eqref{eqn:ref_dis} as follows
\begin{equation} \label{eqn:ref_dis_noise}
	\mbfd_{\lambda(t)}(f_N,g)=W_2^2\brac{\frac{f_N^2}{\innprod{f_N^2}},\frac{g^2+\lambda}{\innprod{g^2+\lambda}}}.
\end{equation}
Here $\lambda=\lambda(t)$ is a given function of $t$ satisfying
\begin{equation*}
	\lambda(t)+g^2(t)>0, \quad t\in[0,t_f].
\end{equation*}
Let $\lambda(t)=2\mu g(t)+\mu^2+\sigma^2$, then we have
\begin{multline*}
	\mbfd_{\lambda}(f_N,g)=W_2^2\brac{\frac{f_N^2}{\innprod{f_N^2}},\frac{g^2+\lambda}{\innprod{g^2+\lambda}}}
		=W_2^2\brac{\frac{g^2+2gr_N+r_N^2}{\innprod{g^2+2gr_N+r_N^2}},\frac{g^2+\lambda}{\innprod{g^2+\lambda}}} \\
		=W_2^2\brac{\frac{(g^2+\lambda)+(2gr_N+r_N^2-\lambda)}{\innprod{(g^2+\lambda)+(2gr_N+r_N^2-\lambda)}},\frac{g^2+\lambda}{\innprod{g^2+\lambda}}}.
\end{multline*}
Applying Theorem \ref{thm:noise1}, we obtain
\begin{equation*}
	\mbe \mbfd_{\lambda}(f_N,g)=O(\frac{1}{N}).
\end{equation*}

\begin{remark}
	In many practical problems, $\mbe r_j=\mu=0$. In such a case, $\lambda=\sigma^2$ is a constant independent of the variable $t$.
\end{remark}

\begin{example} \label{exam:noise_uni}
	In this example, we investigate the influence caused by the uniform distribution $r_j\sim U[-0.1,0.1],\;j=1,2,\cdots,N$. It follows that
	\begin{equation*}
		\mbe r_j=0, \quad \mbd r_j=\frac{1}{300}.
	\end{equation*}
	Let $g(t)$ be the Ricket wavelet $R(t-2.5)$ with $A=1$ and $f_0=2$Hz. The signal $f_N(t)$ and the noise function $r_N(t)$ are given in \eqref{eqn:fNg}-\eqref{eqn:rNt}. The time interval is $[0,5]$. According to the above discussion, $\lambda_*=\frac{1}{300}$. In Table \ref{tab:noise_uni}, we output the expectation values of the distance $\mbe \mbfd_{\lambda}(f_N,g)$ with respect to the parameter $\lambda$ and the number of time divisions $N$. For each configuration, we repeat $100$ trials to compute the expectation values. For reference, we also output the expectation values of the $L^2$ distance between $f_N(t)$ and $g(t)$,
	\begin{equation*}
		\norm{f_N(t)-g(t)}_2=\brac{\int_0^5\abs{f_N(t)-g(t)}^2\rd t}^{1/2}.
	\end{equation*}
	From the table, we can see that $\mbe \mbfd_{\lambda}(f_N,g)\approx O(\frac{1}{N})$ when $\lambda=\lambda_*$. This agrees with our theoretical discussion. Moreover,  $\mbe \mbfd_{\lambda}(f_N,g)$ decreases as $N$ increases when $\lambda$ is close to $\lambda_*$. On the other hand, the expectation values of the $L^2$ distance remains unchanged.
\end{example}

\begin{table*}
    	\caption{Example \ref{exam:noise_uni}: the expectation values of the distance $\mbe \mbfd_{\lambda}(f_N,g)$ with respect to $\lambda$ and $N$. The last line is the expectation values of the $L^2$ distance between $f_N(t)$ and $g(t)$.} \label{tab:noise_uni}
	\begin{center}\begin{tabular}{c|ccccc} \hline 
		\diagbox[dir=SE]{$\lambda$}{$N$} & $50$ & $100$ & $200$ & $400$ & $800$ \\ \hline
		$0.8\lambda_*$ & $1.02\times10^{-2}$ & $8.43\times10^{-3}$ & $6.36\times10^{-3}$ & $5.56\times10^{-3}$ & $5.09\times10^{-3}$ \\ 
		$0.9\lambda_*$ & $8.65\times10^{-3}$ & $4.80\times10^{-3}$ & $2.79\times10^{-3}$ & $2.14\times10^{-3}$ & $1.64\times10^{-3}$ \\ 
		$\lambda_*$ & $7.42\times10^{-3}$ & $4.10\times10^{-3}$ & $2.09\times10^{-3}$ & $9.90\times10^{-4}$ & $5.34\times10^{-4}$ \\ 
		$1.1\lambda_*$ & $6.89\times10^{-3}$ & $4.78\times10^{-3}$ & $3.03\times10^{-3}$ & $1.99\times10^{-3}$ & $1.63\times10^{-3}$ \\
		$1.2\lambda_*$ & $1.03\times10^{-2}$ & $7.01\times10^{-3}$ & $5.31\times10^{-3}$ & $4.79\times10^{-3}$ & $4.04\times10^{-3}$ \\ \hline
		$\mbe\norm{f_N(t)-g(t)}_2$ & $1.70\times10^{-2}$ & $1.66\times10^{-2}$ & $1.66\times10^{-2}$ & $1.67\times10^{-2}$ & $1.67\times10^{-2}$ \\ \hline
    	\end{tabular}\end{center}
\end{table*}

\begin{example} \label{exam:noise_norm}
	In this example, we investigate the influence caused by the normal distribution $r_j\sim N[0,0.1^2],\;j=1,2,\cdots,N$. It follows that
	\begin{equation*}
		\mbe r_j=0, \quad \mbd r_j=0.01.
	\end{equation*}
	Let $g(t)$ be the Ricket wavelet $R(t-2.5)$ with $A=1$ and $f_0=2$Hz. The signal $f_N(t)$ and the noise function $r_N(t)$ are given in \eqref{eqn:fNg}-\eqref{eqn:rNt}. The time interval is $[0,5]$. According to the above discussion, $\lambda_*=\frac{1}{100}$. In Table \ref{tab:noise_norm}, we output the expectation values of the distance $\mbe \mbfd_{\lambda}(f_N,g)$ with respect to the parameter $\lambda$ and the number of time divisions $N$. For each configuration, we repeat $100$ trials to compute the expectation values. For reference, we also output the expectation values of the $L^2$ distance between $f_N(t)$ and $g(t)$. From the table, we can draw the same conclusion as in Example \ref{exam:noise_uni}.
\end{example}

\begin{table*}
    	\caption{Example \ref{exam:noise_norm}: the expectation values of the distance $\mbe \mbfd_{\lambda}(f_N,g)$ with respect to $\lambda$ and $N$. The last line is the expectation values of the $L^2$ distance between $f_N(t)$ and $g(t)$.} \label{tab:noise_norm}
	\begin{center}\begin{tabular}{c|ccccc} \hline 
		\diagbox[dir=SE]{$\lambda$}{$N$} & $50$ & $100$ & $200$ & $400$ & $800$ \\ \hline
		$0.8\lambda_*$ & $4.77\times10^{-2}$ & $2.73\times10^{-2}$ & $1.69\times10^{-2}$ & $1.53\times10^{-2}$ & $1.06\times10^{-2}$ \\ 
		$0.9\lambda_*$ & $4.45\times10^{-2}$ & $2.30\times10^{-2}$ & $1.40\times10^{-2}$ & $7.28\times10^{-3}$ & $4.97\times10^{-3}$ \\ 
		$\lambda_*$ & $3.74\times10^{-2}$ & $2.01\times10^{-2}$ & $1.26\times10^{-2}$ & $6.30\times10^{-3}$ & $3.00\times10^{-3}$ \\ 
		$1.1\lambda_*$ & $3.42\times10^{-2}$ & $2.40\times10^{-2}$ & $1.33\times10^{-2}$ & $8.54\times10^{-3}$ & $4.28\times10^{-3}$ \\
		$1.2\lambda_*$ & $3.86\times10^{-2}$ & $2.52\times10^{-2}$ & $1.56\times10^{-2}$ & $9.73\times10^{-3}$ & $9.11\times10^{-3}$ \\ \hline
		$\mbe\norm{f_N(t)-g(t)}_2$ & $4.89\times10^{-2}$ & $4.99\times10^{-2}$ & $5.07\times10^{-2}$ & $5.04\times10^{-2}$ & $4.99\times10^{-2}$ \\ \hline
    	\end{tabular}\end{center}
\end{table*}

The aforementioned discussions and experiments point out that the parameter $\lambda$ should be specified to reduce the impact of noise. This requires us to estimate the mean and variance of the noise, which will cost some extra efforts. Fortunately, there are many statistical methods to estimate these values, e.g. Independent Component Analysis \cite{HyKaOj:01}. Moreover, the results of the numerical experiments show that the estimation does not need to be particularly accurate.

\section{The application to earthquake location} \label{sec:ael}
Up to now, we have proposed the reformulated distance \eqref{eqn:ref_dis} to measure two earthquake signals  and studied its properties. Next, we would like to apply this distance to determine the real earthquake hypocenter $\bxi_T$ and the origin time $\tau_T$. Its mathematical formulation is written as follows
\begin{equation} \label{eqn:inv_prob}
	(\bxi_T,\tau_T)=\argmin_{\bxi,\tau}\sum_r\chi_r(\bxi,\tau),
\end{equation}
where the misfit function at the $r-$th receiver $\chi_r(\bxi,\tau)$ is defined by
\begin{equation} \label{eqn:mist_fun}
	\chi_r(\bxi,\tau)=\mbfd(d_r(t),s(\boeta_r,t)).
\end{equation}
The real earthquake signal $d_r(t)$ and the synthetic earthquake signal $s(\bx,t)$ can be considered as the solution
\begin{equation} \label{eqn:rel_syn}
	d_r(t)=u(\boeta_r,t;\bxi_T,\tau_T), \quad s(\bx,t)=u(\bx,t;\bxi,\tau),
\end{equation}
of the acoustic wave equation initial-boundary value problem
\begin{align} 
	& \frac{\partial^2 u(\bx,t;\bxi,\tau)}{\partial t^2}=\nabla\cdot\brac{c^2(\bx)\nabla u(\bx,t;\bxi,\tau)}
		+R(t-\tau)\delta(\bx-\bxi), \quad \bx,\bxi\in\Omega, \label{eqn:wave} \\
	& u(\bx,0;\bxi,\tau)=\partial_t u(\bx,0;\bxi,\tau)=0, \quad \bx\in\Omega, \label{ic:wave} \\
	& \bn\cdot \brac{c^2(\bx)\nabla u(\bx,t;\bxi,\tau)}=0, \quad \bx\in\partial \Omega. \label{bc:wave}
\end{align}
In the above equations, $c(\bx)$ denotes the wave speed and $\boeta_r$ is the location of the $r-$th receiver. The computational domain $\Omega$ is a subset of the $d-$dimensional real Euclidean space $\br^d$ and $\bn$ is the outward unit normal vector to the domain $\Omega$. In this study, the seismic rupture is modeled by the point source $\delta(\bx-\bxi)$ since its scale is much smaller compared to the scale of seismic wave propagated \cite{AkRi:80, Ma:15}. We also remark that the reflection boundary condition \eqref{bc:wave} is used to simplify the model. There is no essential difficulty to consider other boundary conditions, e.g. the perfectly matched later absorbing boundary condition \cite{KoTr:03}.

\begin{remark} \label{rem:noise}
	In practice, the real signal is superimposed with noise
	\begin{equation*}
		d_r(t)=u(\boeta_r,t;\bxi_T,\tau_T)+r_N(t).
	\end{equation*}
	Thus, we prefer to use the distance given in \eqref{eqn:ref_dis_noise} to define the misfit function 
	\begin{equation} \label{eqn:rede_mis}
		\chi_r(\bxi,\tau)=\mbfd_{\lambda}(d_r(t),s(\boeta_r,t)).
	\end{equation}
\end{remark}

\subsection{The adjoint method}
The perturbation of earthquake hypocenter $\delta \bxi\ll1$ and origin time $\delta \tau\ll1$ would generate the perturbation of wave function
\begin{equation} \label{rel:pert_s}
	\delta s(\bx,t)=u(\bx,t;\bxi+\delta \bxi,\tau+\delta \tau)-u(\bx,t;\bxi,\tau).
\end{equation}
According to \eqref{eqn:wave}-\eqref{bc:wave}, $\delta s(\bx,t)$ satisfies the acoustic wave equation
\begin{align} 
	& \frac{\partial^2 \delta s(\bx,t)}{\partial t^2}=\nabla\cdot\brac{c^2(\bx)\nabla \delta s(\bx,t)} \label{eqn:pert_s} \\
	& \quad\quad +R(t-(\tau+\delta \tau))\delta(\bx-(\bxi+\delta \bxi))-R(t-\tau)\delta(\bx-\bxi), \quad \bx,\bxi\in\Omega,  \nn \\
	& \delta s(\bx,0)=\partial_t \delta s(\bx,0)=0, \quad \bx\in\Omega, \label{ic:pert_s} \\
	& \bn\cdot \brac{c^2(\bx)\nabla \delta s(\bx,0)}=0, \quad \bx\in\partial \Omega. \label{bc:pert_s}
\end{align}
Multiply an arbitrary test funciton $w_r(\bx,t)$ on equation \eqref{eqn:pert_s}, integrate it on $\Omega\times[0,t_f]$  and use integration by parts, we obtain
{\small\begin{multline} \label{eqn:integration}
	\int_0^{t_f}\int_{\Omega}\frac{\partial^2 w_r}{\partial t^2}\delta s\rd \bx\rd t
		-\int_{\Omega}\left.\frac{\partial w_r}{\partial t}\delta s\right|_{t=t_f}\rd \bx
		+\int_{\Omega}\left.w_r\frac{\partial \delta s}{\partial t}\right|_{t=t_f}\rd \bx \\
	=\int_0^{t_f}\int_{\Omega}\delta s\nabla\cdot(c^2\nabla w_r)\rd \bx \rd t
		-\int_0^{t_f}\int_{\partial \Omega}\bn\cdot (c^2\nabla w_r)\delta s\rd \zeta \rd t
                 +\int_0^{t_f} R(t-(\tau+\delta \tau))w_r(\bxi+\delta \bxi,t)-R(t-\tau)w_r(\bxi,t)\rd t \\
	\approx\int_0^{t_f} \int_{\Omega}\delta s\nabla\cdot(c^2\nabla w_r)\rd \bx \rd t
                  -\int_0^{t_f} \int_{\partial \Omega}\bn\cdot (c^2\nabla w_r)\delta s\rd \zeta \rd t
                  +\int_0^{t_f} R(t-\tau)\nabla w_r(\bxi,t)\cdot \delta \bxi-R'(t-\tau)w_r(\bxi,t)\delta \tau\rd t.
\end{multline}}
In the last step, the Taylor expansion is used and higher order terms are ignored.

On the other hand, the misfit function \eqref{eqn:mist_fun} also generates the perturbation with respect to $\delta s(\bx,t)$, assume that $\norm{\delta s(\bx,t)}\ll 1$, taking into account of \eqref{eqn:frechet_gradient}, we have
\begin{multline} \label{fun:deltachi}
	\delta \chi_r=\chi_r(\bxi+\delta \bxi,\tau+\delta \tau)-\chi_r(\bxi,\tau) \\
		\approx\int_0^{t_f}4\brac{A(t)-B}s(\boeta_r,t)\delta s(\boeta_r,t)\rd t \\
		=\int_0^{t_f}\int_{\Omega}4\brac{A(t)-B}s(\boeta_r,t)\delta s(\bx,t)\delta(\bx-\boeta_r)\rd \bx\rd t.
\end{multline}
where ``$\approx$'' is obtained by ignoring high order terms of $\delta s(\bx,t)$.

Let $w_r(\bx,t)$ satisfies the adjoint equation
\begin{align} 
	& \frac{\partial^2 w_r(\bx,t)}{\partial t^2}=\nabla\cdot\brac{c^2(\bx)\nabla w_r(\bx,t)} 
		+4\brac{A(t)-B}s(\boeta_r,t)\delta(\bx-\boeta_r), \quad \bx,\bxi\in\Omega, \label{eqn:adjoint} \\
	& w_r(\bx,t_f)=\frac{\partial w_r(\bx,t_f)}{\partial t}=0, \quad \bx\in\Omega, \label{ic:adjoint} \\
	& \bn\cdot\brac{c^2(\bx)\nabla w_r(\bx,t)}=0, \quad \bx\in\partial \Omega. \label{bc:adjoint}
\end{align}
Thus, the relation between $\delta \chi_r$ and $\delta \bxi, \;\delta \tau$ can be obtained by adding \eqref{eqn:integration} to \eqref{fun:deltachi}
\begin{equation} \label{eqn:delchi_rel}
	\delta \chi_r=K_r^{\bxi}\cdot \delta\bxi+ K_r^{\tau}\delta \tau,
\end{equation}
in which the sensitivity kernel for the hypocenter $\bxi$ and the origin time $\tau$ is
\begin{eqnarray}
	K_r^{\bxi} &=& \int_0^{t_f} R(t-\tau)\nabla w_r(\bxi,t)\rd t, \\
	K_r^{\tau} &=& -\int_0^{t_f} R'(t-\tau)w_r(\bxi,t)\rd t.
\end{eqnarray}

\subsection{The LMF method}
According to \eqref{eqn:qwm_sol}, \eqref{eqn:ref_dis} and \eqref{eqn:inv_prob}-\eqref{bc:pert_s}, the mathematical formulation of the earthquake location is a least square optimization problem. Therefore, we can consider some special methods to improve the convergence. Throungh a large number of numerical tests, we found that the LMF method \cite{Fl:91, Le:44, MaNiTi:04, Ma:63} works very well. In the following, we briefly review the basic idea of the algorithm.

In order to be consistent with the literatures of optimization theory, all the symbols and notations in this subsection is independent from the other part of the paper. The general form of the least-square problem can be written as
\begin{equation} \label{eqn:lsp}
	\min f(\bx)=\half \sum_{i=1}^mr_i^2(\bx), \quad \bx\in\br^n,\;m\ge n,
\end{equation}
where the residual function
\begin{equation*}
	\vr(\bx)=(r_1(\bx),r_2(\bx),\cdots,r_m(\bx))^T\in\br^m.
\end{equation*}
The gradient of $f$ is
\begin{equation}
	\nabla f(\bx)=\sum_{i=1}^mr_i(\bx)\nabla r_i(\bx)=J(\bx)^T\vr(\bx),
\end{equation}
in which $J(\bx)$ is the Jacobi matrix of $\vr(\bx)$
\begin{equation*}
	J(\bx)=\brac{\nabla r_1, \nabla r_2, \cdots, \nabla r_m}^T\in\br^{m\times n}.
\end{equation*}

The key ingredient of LMF method is
\begin{equation} \label{eqn:LMF}
	\brac{J_k^TJ_k+\nu_kI}\bd_k=-J_k^T\vr_k, \quad \bx_{k+1}=\bx_k+\bd_k,
\end{equation}
here $I$ is the identity matrix and $\nu_k\ge 0$ is a parameter in each iteration step. It is introduced to improve the convergence and efficiency. To adjust this parameter, we define
\begin{equation} \label{eqn:LMF_ind}
	\gamma_k=\frac{f(\bx_k)-f(\bx_k+\bd_k)}{q_k(0)-q_k(\bd_k)},
\end{equation}
with
\begin{equation}
	q_k(\bd)=\half\brac{J_k\bd+\vr_k}^T\brac{J_k\bd+\vr_k}.
\end{equation}
Now the detailed implementation of the LMF method is summarized below.

\begin{algorithm}[The LMF method ] \label{alg:LMF_nf} $ $
\begin{enumerate}[1.]
	\item Set tolerance value $\veps=0.01$, the break-off step $K=20$ and $\mu=2$. Let $k=0$ and give the initial value $\bx_0$. Thus, we have the initial adjustable parameter $\nu_0=10^{-6}\times\max\abs{\textrm{diag}(J_0^TJ_0)}$.
	
	\item For $\bx_k$ and $\nu_k$, solve the equation \eqref{eqn:LMF} to obtain $\bd_k$. And we can calculate $\gamma_k$ using equation \eqref{eqn:LMF_ind}.
	
	\item If $\gamma_k>0$, let $\bx_{k+1}=\bx_k+\bd_k,\;\nu_{k+1}=\nu_k\times\max\{\frac{1}{3},1-(2\gamma_k-1)^3\},\;\mu=2$ and $k=k+1$. If $\norm{f(\bx_k)}<\veps$, output $\bx_k$ and stop. 
	
	\item If $\gamma_k<0$, let $\nu_k=\mu\nu_k,\;\mu=2\mu$.
	
	\item If $k>K$, output the error message:``The iteration doesn't converges.'' and stop. Otherwise go to step 2 for another iteration. 
\end{enumerate}
\end{algorithm}

In the above algorithm, we require that the objective function always decrease. Otherwise, we will decrease the radius of the trust region. For the noise-free situation, the idea works well since the objective function has nice convexity in a large area. However, due to the influence of data noise, the optimization objective function will appear some local minimum. Therefore, it may be not suitable to require the objective function decrease during the whole iteration. Thus, we modify the LMF method as follows.

\begin{algorithm}[The modified LMF method for noise data] \label{alg:LMF_nd} $ $
\begin{enumerate}[1.]
	\item Set tolerance value $\veps=0.01$, the break-off step $K=20$ and $\mu=2$. Let $k=0$ and give the initial value $\bx_0$. Thus, we have the adjustable parameter $\nu_0=10^{-3}$ and its upper limit $\eta=10^{-3}$.
	
	\item For $\bx_k$ and $\nu_k$, solve the equation \eqref{eqn:LMF} to obtain $\bd_k$. And we can calculate $\gamma_k$ using equation \eqref{eqn:LMF_ind}.
	
	\item If $\gamma_k>0$ or $\nu_k\ge \eta$, let $\bx_{k+1}=\bx_k+\bd_k,\;\nu_{k+1}=\min\{\eta,\nu_k\times\max\{\frac{1}{3},1-(2\gamma_k-1)^3\}\},\;\mu=2$ and $k=k+1$. If $\norm{f(\bx_k)}<\veps$, output $\bx_k$ and stop. 
	
	\item If $\gamma_k<0$ and $\nu_k<\eta$, let $\nu_k=\min\{\eta,\mu\nu_k\},\;\mu=2\mu$.
	
	\item If $k>K$, output the error message:``The iteration doesn't converges.'' and stop. Otherwise go to step 2 for another iteration. 
\end{enumerate}
\end{algorithm}

\begin{remark}
	The criteria $\norm{f(\bx_k)}<\veps$ in the above algorithm may not be reached because of the data noise. For this problem, we have two options. One is to change $\veps$, but this could be a little tricky to choose an appropriate parameter $\veps$. The other is to wait until the iteration ends. Then, we can choose the numerical solutions corresponding to the smallest value of the objective function
	\begin{equation*}
		k_*=\argmin_{k}\norm{f(\bx_k)}.
	\end{equation*}
\end{remark}

\section{Numerical Experiments} \label{sec:num}
In this section, two examples are presented to demonstrate the validity of the optimal transport model. We first study the convexity of the optimization objective function with respect to the initial earthquake hypocenter. The convergency and the efficiency of Algorithm \ref{alg:LMF_nf} will be discussed for the situation when the seismic signals are noise-free. For the situation of data noise, the convergency of Algorithm \ref{alg:LMF_nd} are numerically investigated.

To solve the acoustic wave equation \eqref{eqn:wave}, we simply use the finite difference schemes \cite{Da:86, LiYaWuMa:17, YaEnSuFr:16}. Inside the earth, the perfectly matched layer boundary condition \cite{KoTr:03} is applied to absorb the outgoing waves. On the other hand, the reflection boundary condition \eqref{bc:wave} is used to model the free surface of the earth. To discretize the delta function for the point source $\delta(\bx-\bxi)$, we borrow the idea from \cite{We:08}. It writes
\begin{equation*}
	\delta(x)=\left\{\begin{array}{ll}
		\frac{1}{h}\left(1-\frac{5}{4}\abs{\frac{x}{h}}^2-\frac{35}{12}\abs{\frac{x}{h}}^3
			+\frac{21}{4}\abs{\frac{x}{h}}^4-\frac{25}{12}\abs{\frac{x}{h}}^5\right), & \abs{x}\le h, \\
		\frac{1}{h}\left(-4+\frac{75}{4}\abs{\frac{x}{h}}-\frac{245}{8}\abs{\frac{x}{h}}^2+\frac{545}{24}\abs{\frac{x}{h}}^3
			-\frac{63}{8}\abs{\frac{x}{h}}^4+\frac{25}{24}\abs{\frac{x}{h}}^5\right), & h<\abs{x}\le 2h, \\
		\frac{1}{h}\left(18-\frac{153}{4}\abs{\frac{x}{h}}+\frac{255}{8}\abs{\frac{x}{h}}^2-\frac{313}{24}\abs{\frac{x}{h}}^3
			+\frac{21}{8}\abs{\frac{x}{h}}^4-\frac{5}{24}\abs{\frac{x}{h}}^5\right), & 2h<\abs{x}\le 3h, \\
		0, & \abs{x}>3h.
	\end{array}\right. 
\end{equation*}
Here $h$ is a numerical parameter which is related to the mesh size.

\subsection{The two-layer model} \label{subsec:twolayer}
\begin{figure} 
	\centering
	\includegraphics[width=0.70\textwidth, height=0.18\textheight]{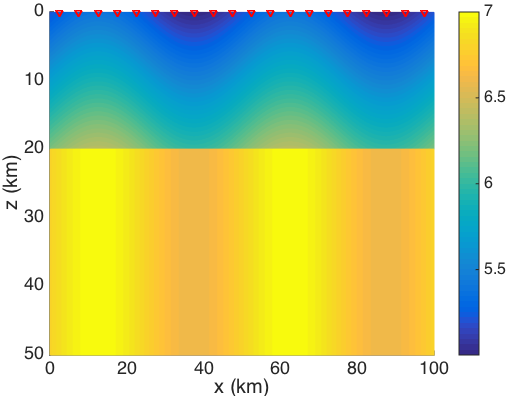}
	\caption{Illustration of two-layer model. The read triangles indicate the receivers.} \label{fig:exam41_vel}
\end{figure}

Consider the two-layer model in the bounded domain $\Omega=[0,100\,km]\times[0,50\,km]$, the wave speed is
\begin{equation*}
	c(x,z)=\left\{\begin{array}{ll}
		5.2+0.05z+0.2\sin\frac{\pi x}{25}, & 0\,km\le z\le 20\,km, \\
		6.8+0.2\sin\frac{\pi x}{25}, & z>20\,km.
	\end{array}\right.
\end{equation*}
The unit is `km/s'. The computational time interval $I=[0,35\,s]$. The dominant frequency of the earthquakes is $f_0=2\,Hz$. There are $20$ equidistant receivers on the surface
\begin{equation*}
	\boeta_r=(x_r,z_r)=(5r-2.5\,km,0), \quad r=1,2,\cdots,20,
\end{equation*}
see Figure \ref{fig:exam41_vel} for illustration.

\begin{figure} 
	\centering
	\includegraphics[width=0.90\textwidth, height=0.32\textheight]{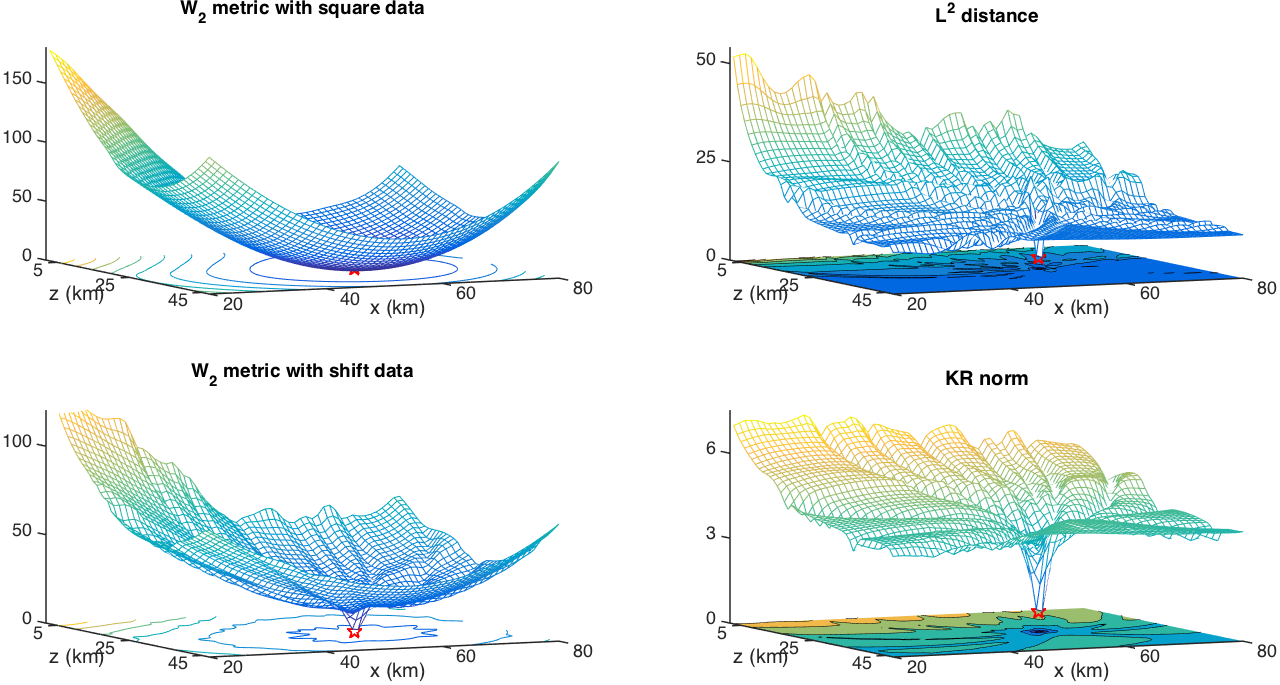}
	\caption{The two-layer model, the cross-section of the optimization objective function with respect to different measures, case (i). The red pentagram denotes the real earthquake hypocenter.} \label{fig:exam41_illu_deep_convex}
\end{figure}

\begin{figure} 
	\centering
	\includegraphics[width=0.90\textwidth, height=0.32\textheight]{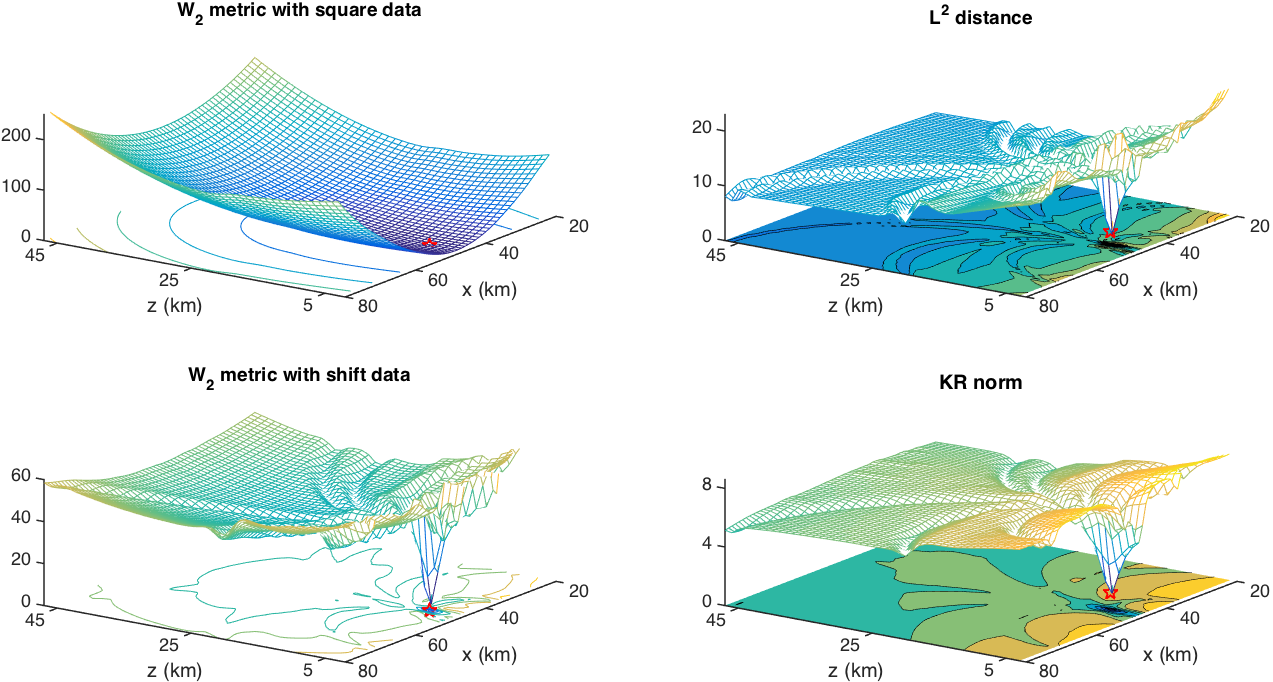}
	\caption{The two-layer model, the cross-section of the optimization objective function with respect to different measures, case (ii). The red pentagram denotes the real earthquake hypocenter.} \label{fig:exam41_illu_shallow_convex}
\end{figure}

First, we output the cross-section of the optimization objective function
\begin{equation*}
	\Psi(\bxi)=\sum_r\chi_r(\bxi,\tau_T),
\end{equation*}
in which $\chi_r(\bxi,\tau)$ is defined in \eqref{eqn:mist_fun}. Here, the distance between the real signal $d_r(t)$ and the synthetic signal $s(\boeta_r,t)$ is measured by the quadratic Wasserstein metric of the normalized square signals (QWN$_2$) 
\begin{equation*}
	W_2^2\brac{\frac{d_r^2(t)}{\innprod{d^2_r(t)}},\frac{s^2(\boeta_r,t)}{\innprod{s^2(\boeta_r,t)}}}.
\end{equation*}
As a comparison, we also output the corresponding objective function under other types of distance as follows:
\begin{itemize}
	\item The relative $L^2$ distance (RLD): 
		\begin{equation*}
			\frac{\int_0^{t_f}\abs{d_r(t)-s(\boeta_r,t)}^2\rd t}{\int_0^{t_f}\abs{d_r(t)}^2\rd t};
		\end{equation*}
	\item The quadratic Wasserstein metric of the normalized shift signals (QWN$_c$) \cite{YaEnSuFr:16}: 
		\begin{equation*}
			W_2^2\brac{\frac{d_r(t)+c}{\innprod{d_r(t)+c}},\frac{s(\boeta_r,t)+c}{\innprod{s(\boeta_r,t)+c}}}, 
		\end{equation*}
		here $c$ is a constant to ensure the positive;
	\item The Kantorovich-Rubinstein norm of the original signals (KRN) \cite{KaRu:58, MeBrMeOuVi:16}:
		\begin{equation*}
			W_1(d_r(t),s(\boeta_r,t))=\max_{\varphi\in \textrm{BLip}_1}\int_0^{t_f}\varphi(t)\brac{d_r(t)-s(\boeta_r,t)}\rd t,
		\end{equation*}
		in which $\textrm{BLip}_1$ is the space of bounded 1-Lipschitz functions, such that
		\begin{equation*}
			(i)\;\forall (t_1,t_2)\in[0,t_f], \; \abs{\varphi(t_1)-\varphi(t_2)}\le \abs{t_1-t_2}, \quad
			(ii)\;\forall t\in[0,t_f],\; \abs{\varphi(t)}\le 1.
		\end{equation*}
\end{itemize}
Here we consider two different earthquake hypocenter, one is below the Moho discontinuity (Figure \ref{fig:exam41_illu_deep_convex}) and another is above the Moho discontinuity (Figure \ref{fig:exam41_illu_shallow_convex}):
\begin{align*}
	& (i) \quad \bxi_T=(57.604\,km,26.726\,km), \quad \tau_T=10.184\,s. \\
	& (ii) \quad \bxi_T=(46.234\,km,7.124\,km), \quad \tau_T=10.782\,s;
\end{align*} 
From these figures, we can observe nice convexity property of the optimization objective function $\Psi(\bxi)$ with respect to earthquake hypocenter $\bxi=(\zeta_x,\zeta_z)$ by QWN$_2$. For other distances, it seems the convexity property is not good enough.

\begin{figure} 
	\centering
	\includegraphics[width=0.90\textwidth, height=0.16\textheight]{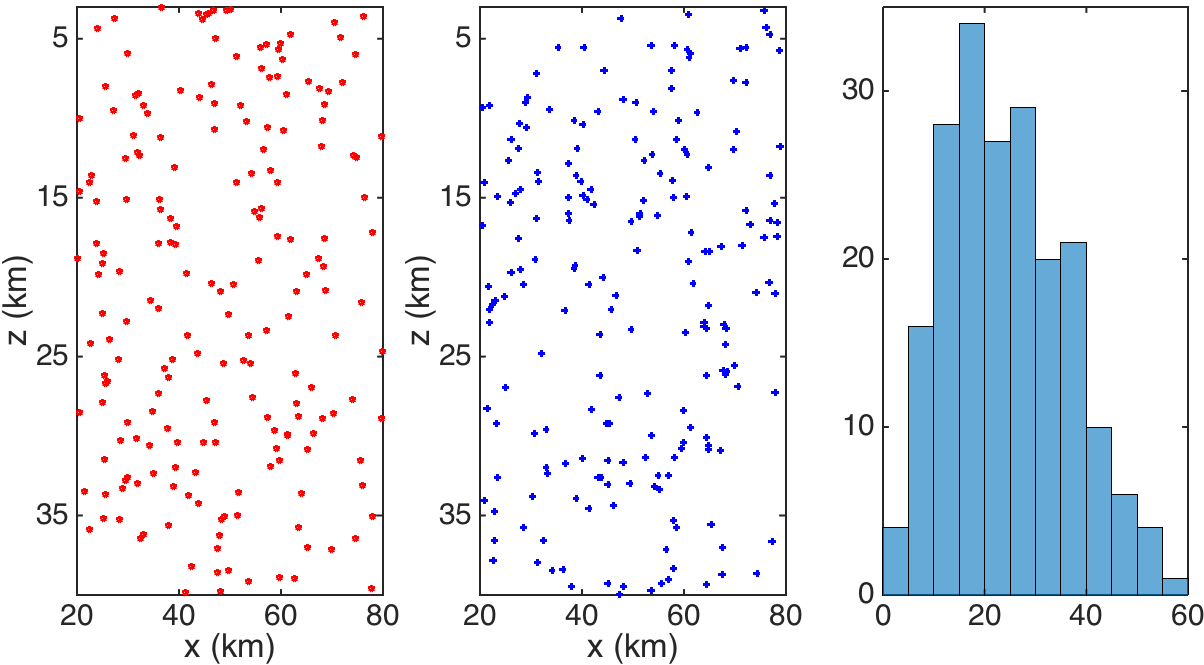}
	\caption{The two-layer model. Left: the spatial distribution of the real earthquake hypocenter $\bxi_T^i$; Middle: the spatial distribution of the initial earthquake hypocenter $\bxi^i$; Right: the distance distribution histogram between the real and the initial earthquake hypocenter $d^i$.} \label{fig:exam41_dis}
\end{figure}

Next, we test the LMF method (Algorithm \ref{alg:LMF_nf}) using $200$ experiments. The real and initial earthquake hypocenter $\bxi_T^i,\;\bxi^i$ are both uniformly distributed over $[20\,km,\;80\,km]\times[3,\;40\,km]$, the real and initial original time $\tau_T^i,\;\tau^i$ are both uniformly distributed over $[7.5\,s,\;12.5\,s]$. Their spatial distribution and the histogram of the distance between the real and the initial hypocenter
\begin{equation*}
	d^i=\norm{\bxi_T^i-\bxi^i}_2,
\end{equation*}
are presented in Figure \ref{fig:exam41_dis}. For all the methods, we randomly select seven receivers for inversion, e.g. $r=4,\,5,\,7,\,9,\,12,\,14,\,18$. In Table \ref{tab:exam41_convergent}, we can see the convergence results for the LMF method, the Gauss-Newton (GN) method and the BFGS method. From which, we can see the LMF method correctly converges in all the tests. But there are $53$ divergent results by the GN method and $10$ error convergence results by the BFGS method. For the convergent cases, we output the mean and standard deviation of iterations for the three methods in Table \ref{tab:exam41_iteration}. It is obvious to see that the LMF method converges faster than the BFGS method. Considering all the above factors, we can conclude that the LMF method is a better choice here.

\begin{table*}
    	\caption{The two-layer model. Convergent results for the LMF method, the GN method and the BFGS method.} \label{tab:exam41_convergent}
	\begin{center}\begin{tabular}{c|ccc|c} \hline
		 & Correct convergence & Divergence & Error convergence & Total \\ \hline
		 LMF & $200$ & $0$ & $0$ & $200$\\ 
		 GN & $147$ & $53$ & $0$ & $200$ \\ 
		 BFGS & $190$ & $0$ & $10$ & $200$ \\ \hline
    	\end{tabular}\end{center}
\end{table*}

\begin{table*}
    	\caption{The two-layer model. Mean and Standard Deviation of iterations for the LMF method, the GN method and the BFGS method.} \label{tab:exam41_iteration}
	\begin{center}\begin{tabular}{c|cc} \hline
		 & Mean of iterations & Standard Deviation of iterations \\ \hline
		 LMF & $5.93$ & $1.90$\\ 
		 GN & $5.59$ & $1.71$ \\ 
		 BFGS & $10.80$ & $2.84$ \\ \hline
    	\end{tabular}\end{center}
\end{table*}

\begin{figure*} 
	\begin{tabular}{c}
		\includegraphics[width=0.95\textwidth, height=0.12\textheight]{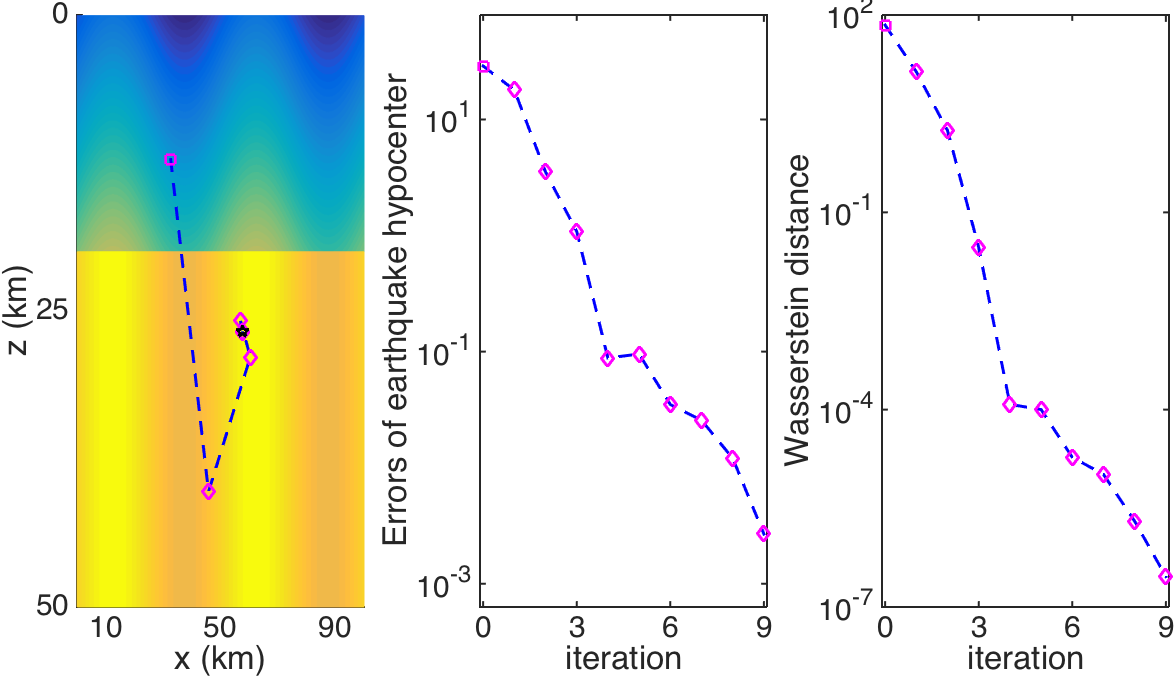} \\
		\includegraphics[width=0.95\textwidth, height=0.12\textheight]{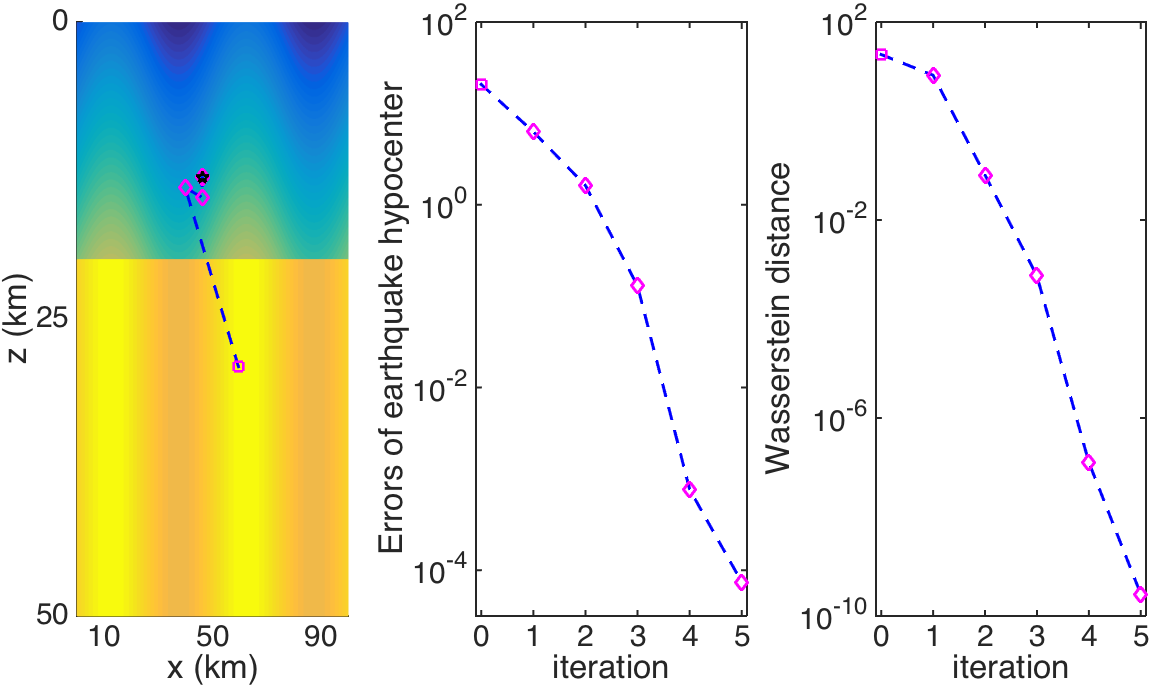}
	\end{tabular}
	\caption{Convergence history of the two-layer model.  Up for case (i), and Down for case (ii). Left: the convergent trajectories; Mid: the absolute errors between the real and computed earthquake hypocenter with respect to iteration steps; Right: the Wasserstein distance between the real and synthetic earthquake signals with respect to iteration steps. The magenta square is the initial hypocenter, the magenta diamond denotes the hypocenter in the iterative process, and the black pentagram is the real hypocenter.} \label{fig:exam41_traj}
\end{figure*}

Then, we output the convergent history of the LMF method (Algorithm \ref{alg:LMF_nf}) by two special examples. Their parameters are selected as follows:
\begin{align*}
	(i) \; \bxi_T=(57.604\,km,\;26.726\,km),\;\tau_T=10.184\,s, \quad \bxi=(32.653\,km,\;12.214\,km),\;\tau=12.108\,s; \\
	(ii) \; \bxi_T=(46.234\,km,\;13.124\,km),\;\tau_T=10.782\,s, \quad \bxi=(59.572\,km,\;29.013\,km),\;\tau=9.908\,s;
\end{align*}
In Figure \ref{fig:exam41_traj}, we can see the convergent trajectories, the absolute errors of the earthquake hypocenter and the Wasserstein distance. These figures show that the method converges to the real earthquake hypocenter very quickly.

\begin{figure*} 
	\begin{tabular}{cccc}
		\includegraphics[width=0.23\textwidth, height=0.12\textheight]{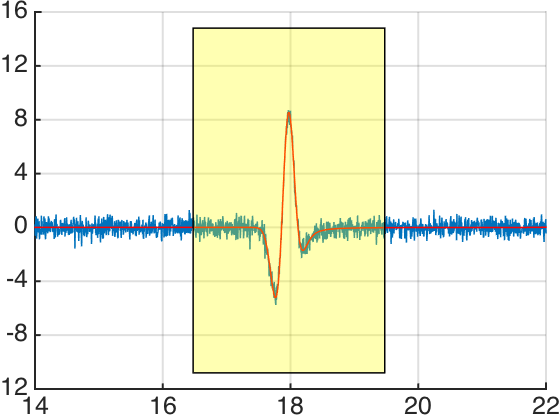} &
		\includegraphics[width=0.23\textwidth, height=0.12\textheight]{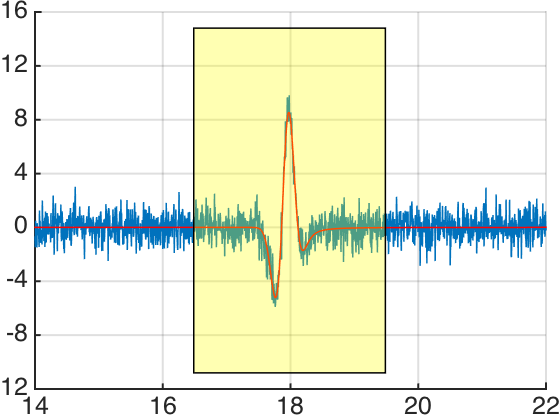} &
		\includegraphics[width=0.23\textwidth, height=0.12\textheight]{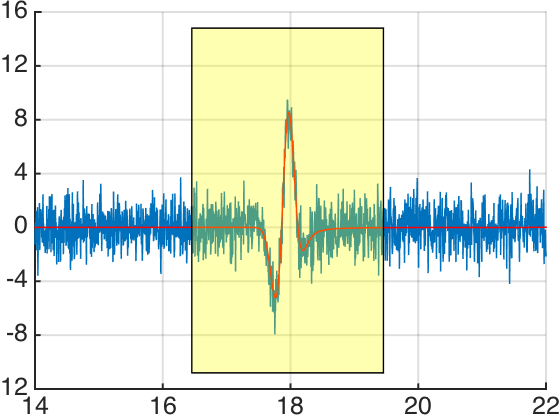} &
		\includegraphics[width=0.23\textwidth, height=0.12\textheight]{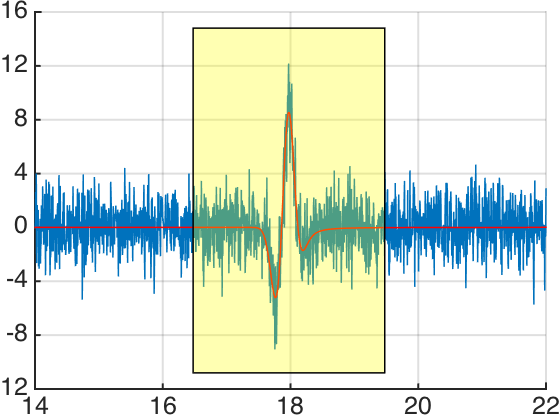} \\
		\includegraphics[width=0.23\textwidth, height=0.12\textheight]{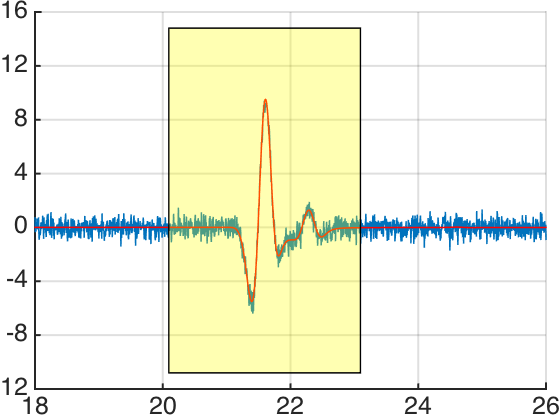} &
		\includegraphics[width=0.23\textwidth, height=0.12\textheight]{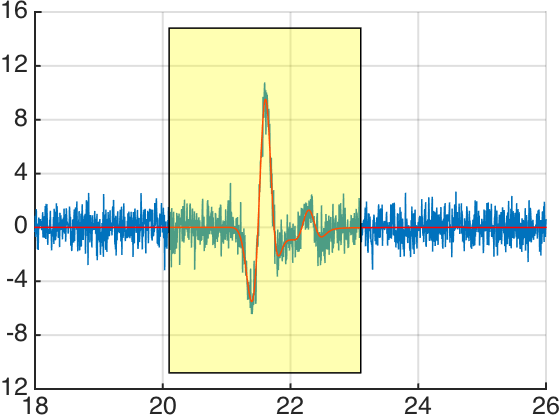} &
		\includegraphics[width=0.23\textwidth, height=0.12\textheight]{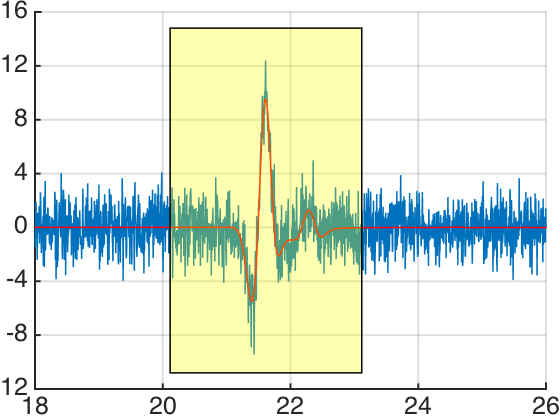} &
		\includegraphics[width=0.23\textwidth, height=0.12\textheight]{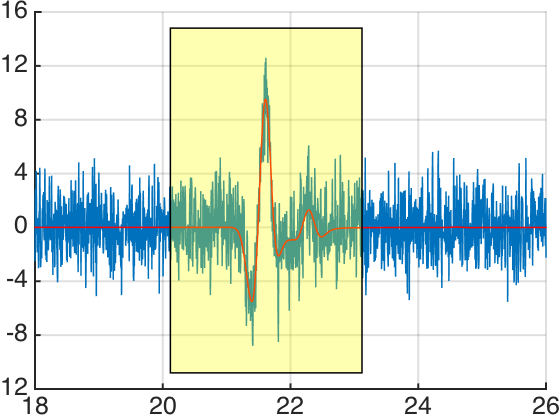}
	\end{tabular}
	\caption{Illustration of signal with noise in the two-layer model. The signal with noise $d_r(t)$ (blue line) and the noise free signal $u(\boeta_r,t;\bxi_T,\tau_T)$ for receiver $r=7$. The horizontal axis is the time $t$. Up: parameters (i); Down: parameters (ii); From left to right, the ratio $R=5\%,\;10\%,\;15\%,\;20\%$ respectively.} \label{fig:exam41_noise_signal}
\end{figure*}

\begin{figure*} 
	\begin{tabular}{c}
		\includegraphics[width=0.95\textwidth, height=0.10\textheight]{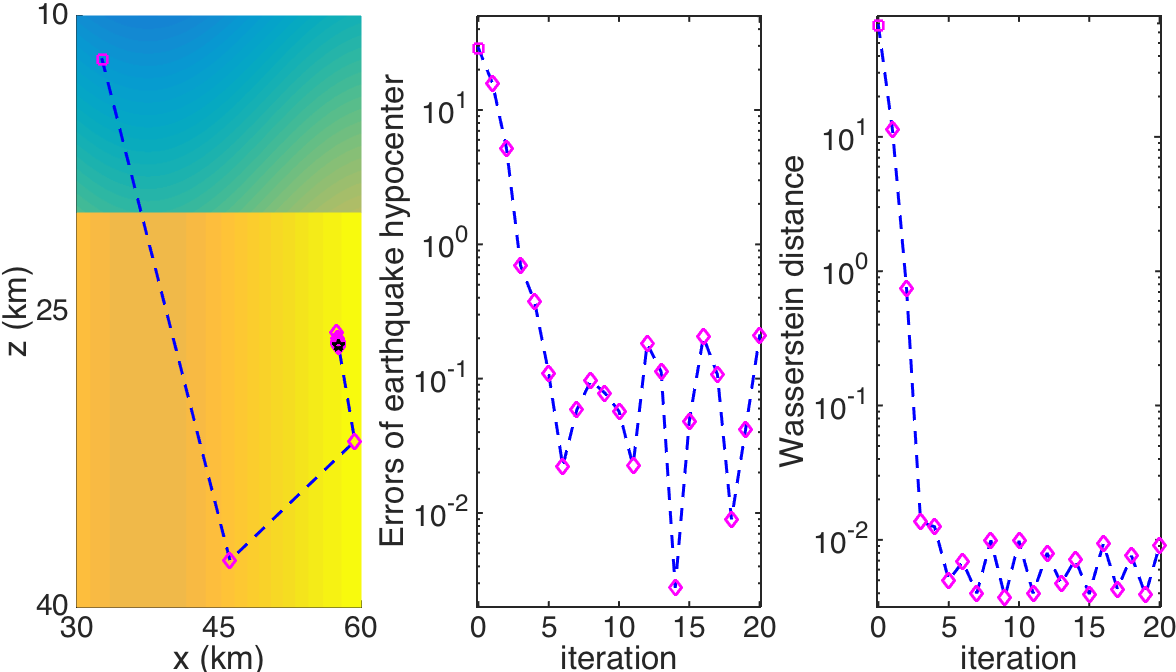} \\
		\includegraphics[width=0.95\textwidth, height=0.10\textheight]{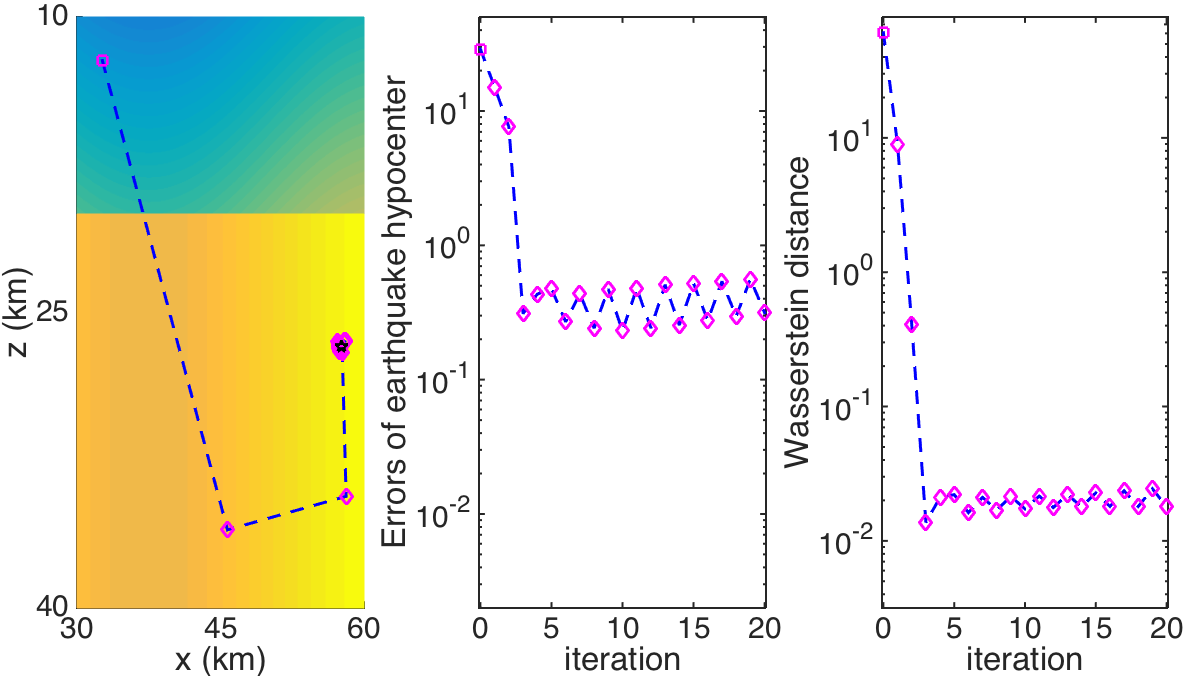} \\
		\includegraphics[width=0.95\textwidth, height=0.10\textheight]{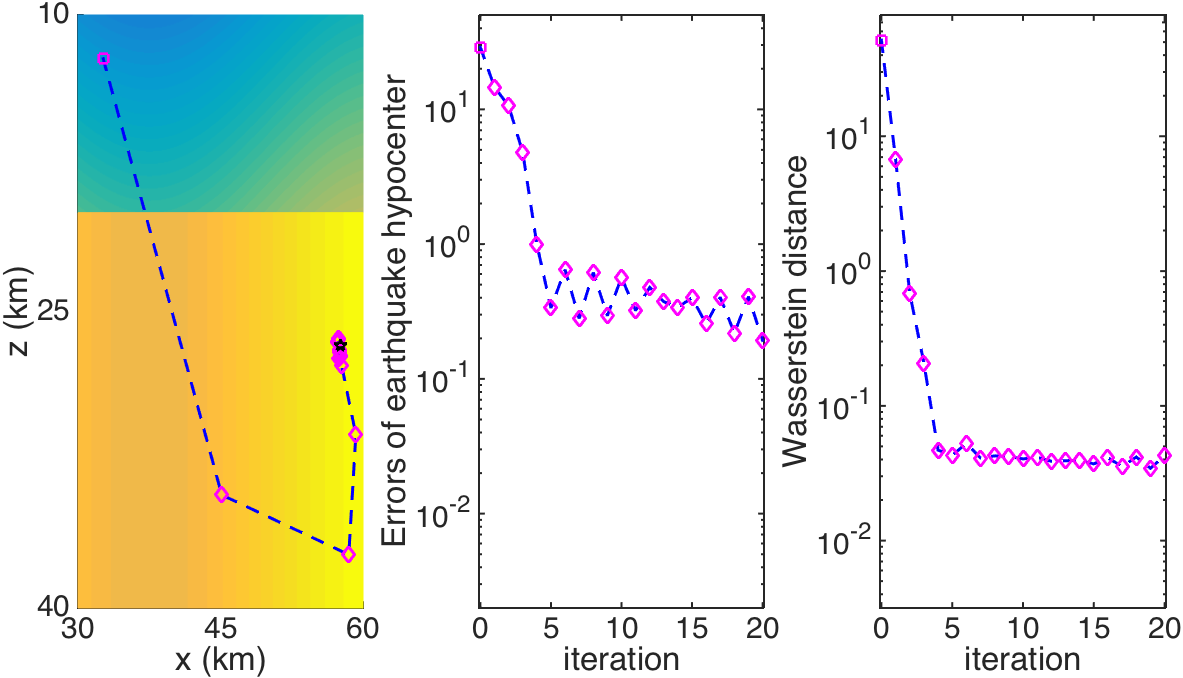} \\
		\includegraphics[width=0.95\textwidth, height=0.10\textheight]{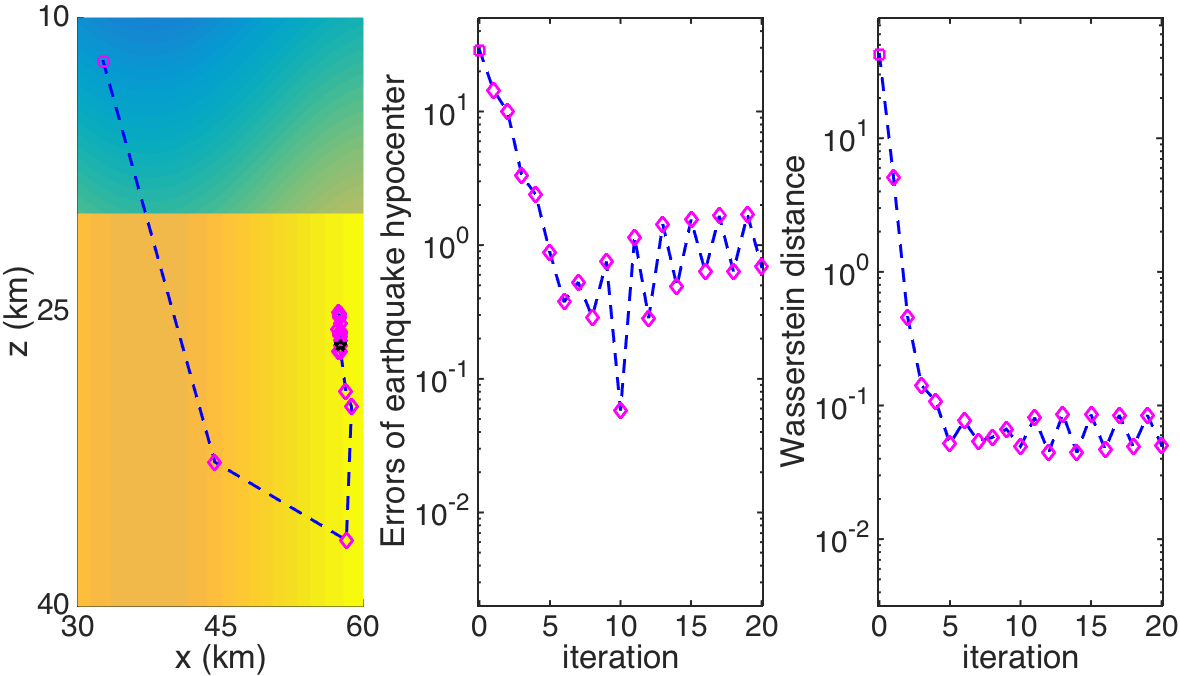}
	\end{tabular}
	\caption{Convergence history of the two-layer model with noise data, case (i). From up to bottom, the the ratio $R=5\%,\;10\%,\;15\%,\;20\%$ respectively. Left: the convergent trajectories; Mid: the absolute errors between the real and computed earthquake hypocenter with respect to iteration steps; Right: the Wasserstein distance between the real and synthetic earthquake signals with respect to iteration steps. The magenta square is the initial hypocenter, the magenta diamond denotes the hypocenter in the iterative process, and the black pentagram is the real hypocenter.} \label{fig:exam41_noise1_traj}
\end{figure*}

\begin{figure*} 
	\begin{tabular}{c}
		\includegraphics[width=0.95\textwidth, height=0.10\textheight]{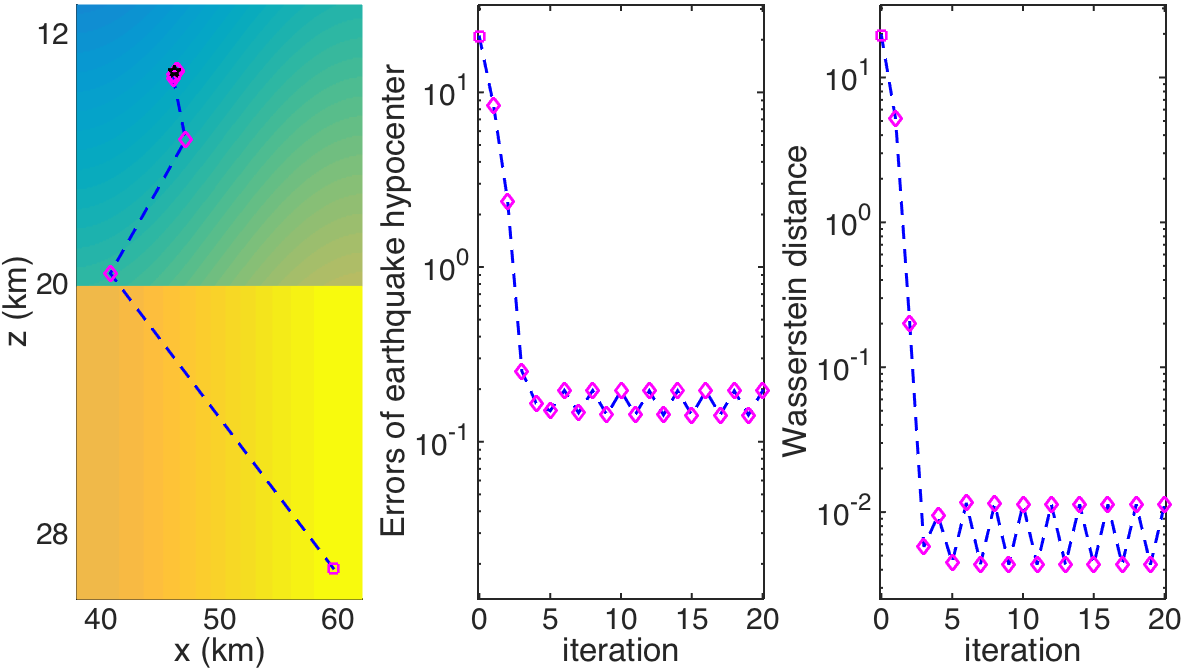} \\
		\includegraphics[width=0.95\textwidth, height=0.10\textheight]{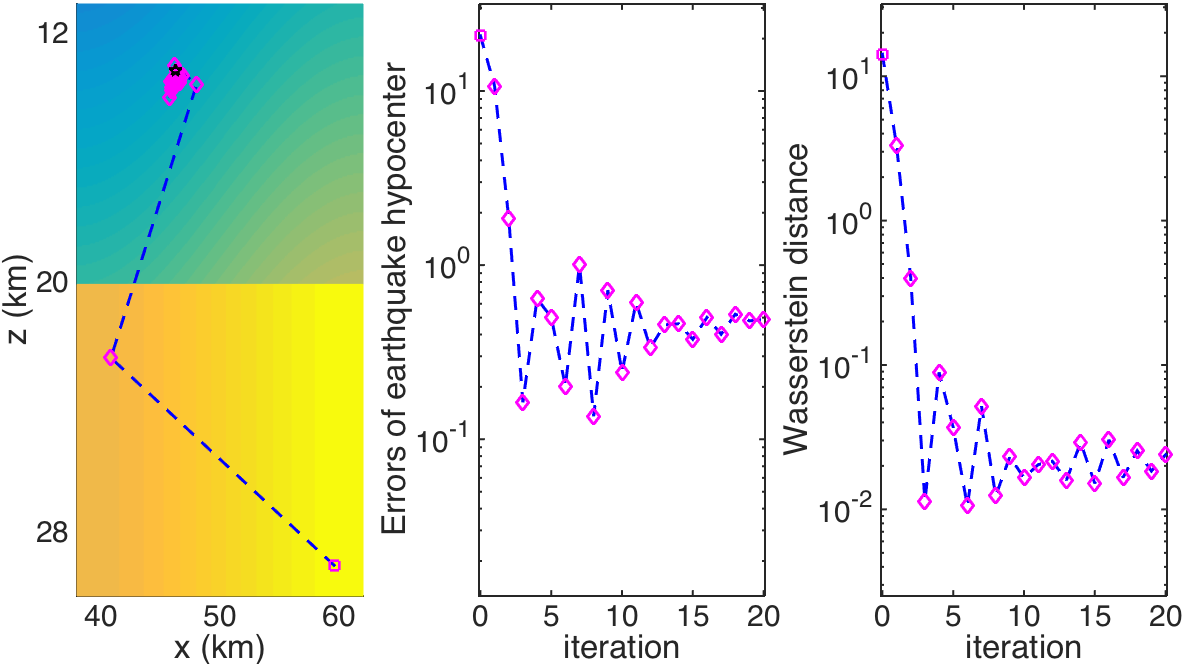} \\
		\includegraphics[width=0.95\textwidth, height=0.10\textheight]{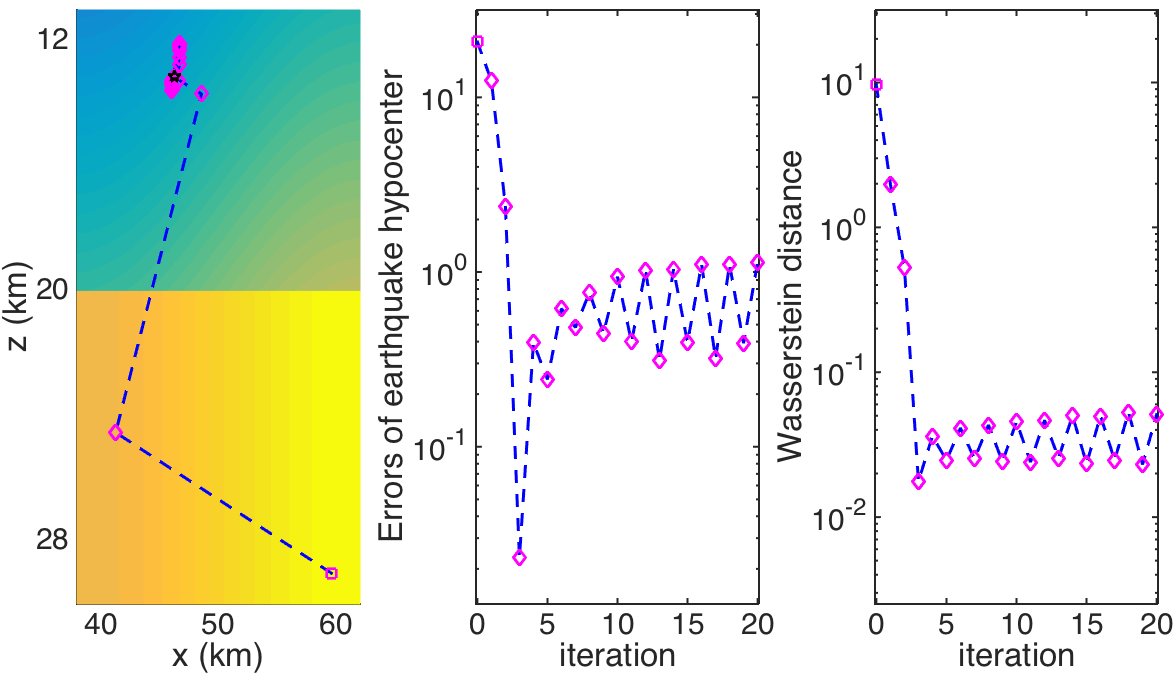} \\
		\includegraphics[width=0.95\textwidth, height=0.10\textheight]{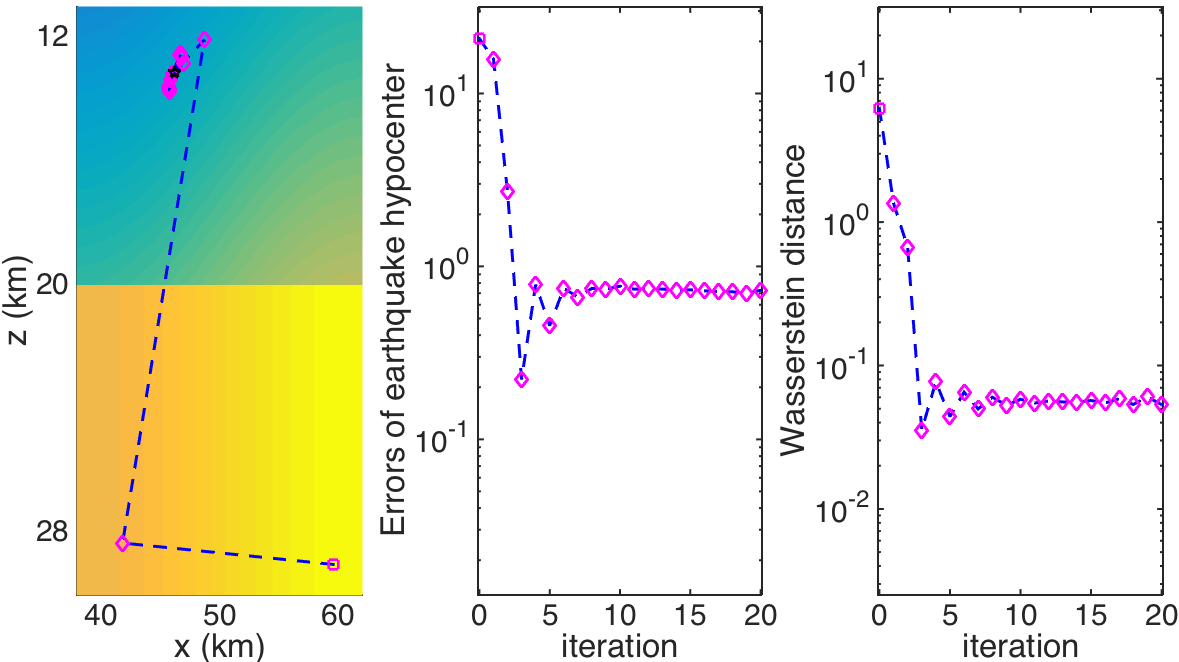}
	\end{tabular}
	\caption{Convergence history of the two-layer model with noise data, case (ii). From up to bottom, the the ratio $R=5\%,\;10\%,\;15\%,\;20\%$ respectively. Left: the convergent trajectories; Mid: the absolute errors between the real and computed earthquake hypocenter with respect to iteration steps; Right: the Wasserstein distance between the real and synthetic earthquake signals with respect to iteration steps. The magenta square is the initial hypocenter, the magenta diamond denotes the hypocenter in the iterative process, and the black pentagram is the real hypocenter.} \label{fig:exam41_noise2_traj}
\end{figure*}

At last, we test the effectiveness of the new method for the data noise. The same parameters (i) and (ii) are selected here. The real earthquake signal can be regarded as
\begin{equation*}
	d_r(t)=u(\boeta_r,t;\bxi_T,\tau_T)+N_r(t).
\end{equation*}
Here $N_r(t)$ is subject to the normal distribution with mean $\mu=0$ and the standard deviation
\begin{equation*}
	\sigma=R \times \max_t\abs{u(\boeta_r,t;\bxi_T,\tau_T)}.
\end{equation*}
The ratio $R$ will be selected as $5\%,\;10\%,\;15\%$ and $20\%$ respectively in the later tests. These signals are illustrated in Fig \ref{fig:exam41_noise_signal}. Obviously, a time window that contains the main part of $u(\boeta_r,t;\bxi_T,\tau_T)$ can be chosen to reduce the impact of noise. As discussed above, we use the formulation \eqref{eqn:rede_mis} in Remark \ref{rem:noise} and the modified LMF method (Algorithm \ref{alg:LMF_nd}) to deal with the noise. The convergent history are output in Figure \ref{fig:exam41_noise1_traj}-\ref{fig:exam41_noise2_traj}.  From these figures, we can see the location errors and the misfit functions oscillate during the iteration. And the iteration step $k_*$ which corresponds to the smallest value of the misfit function does not corresponding to the smallest location error. These are the unavoidable effects of noise. In Table \ref{tab:exam41_noise1}-\ref{tab:exam41_noise2}, we output $k_*$, the corresponding misfit value and the corresponding location error. We can see the location results are still good enough. In addition, we notice that the locations error may be reduced as the noise ratio $R$ increase. This is also caused by the randomness of the noise. Nevertheless, the above results can show a strong adaptability to noise of the new model and method.

\begin{table*}
    	\caption{The two-layer model with noise data, case (i). The smallest misfit value, the corresponding iteration step $k_*$ and the location error.} \label{tab:exam41_noise1}
	\begin{center}\begin{tabular}{c|ccc} \hline
		 $R$ & $k_*$ & The misfit value & The location error (km) \\ \hline
		 $5\%$ & $9$ & $3.74\times 10^{-3}$ & $7.80\times 10^{-2}$ \\ 
		 $10\%$ & $3$ & $1.37\times 10^{-2}$ & $3.10\times 10^{-1}$ \\ 
		 $15\%$ & $10$ & $3.44\times 10^{-2}$ & $4.08\times 10^{-1}$ \\ 
		 $20\%$ & $12$ & $4.45\times 10^{-2}$ & $2.85\times 10^{-1}$ \\ \hline
    	\end{tabular}\end{center}
\end{table*}

\begin{table*}
    	\caption{The two-layer model with noise data, case (ii). The smallest misfit value, the corresponding iteration step $k_*$ and the location error.} \label{tab:exam41_noise2}
	\begin{center}\begin{tabular}{c|ccc} \hline
		 $R$ & $k_*$ & The misfit value & The location error (km) \\ \hline
		 $5\%$ & $15$ & $4.33\times 10^{-3}$ & $1.42\times 10^{-1}$ \\ 
		 $10\%$ & $6$ & $1.06\times 10^{-2}$ & $2.02\times 10^{-1}$ \\ 
		 $15\%$ & $3$ & $1.76\times 10^{-2}$ & $2.31\times 10^{-2}$ \\ 
		 $20\%$ & $3$ & $3.55\times 10^{-2}$ & $2.22\times 10^{-1}$ \\ \hline
    	\end{tabular}\end{center}
\end{table*}

\subsection{The subduction plate model}  \label{subsec:subduction}
\begin{figure} 
	\centering
	\includegraphics[width=0.70\textwidth, height=0.18\textheight]{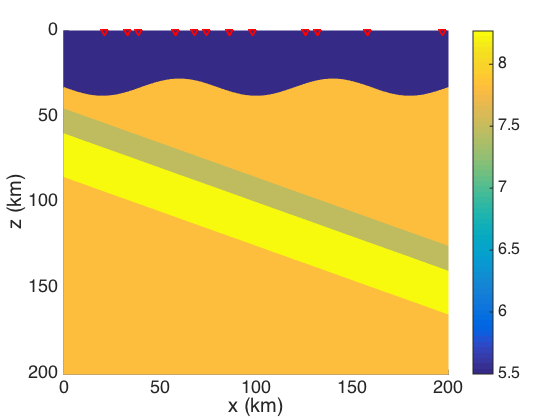}
	\caption{Illustration of the subduction plate model. The read triangles indicate the receivers.} \label{fig:exam42_vel}
\end{figure}

Let's consider a typical seismogenic zone \cite{ToYaLiYaHa:16, ToZhYa:11}. It consists of the crust, the mantle and the undulating Moho discontinuity. And there is a subduction zone with a thin low velocity layer atop a fast velocity layer in the mantle. The earthquake may occur in any of these areas. Taking into account the complex velocity structure, it is much difficult to locate earthquake. In the simulating domain $\Omega=[0,\,200km]\times[0,\,200km]$, the wave speed is 
\begin{equation*}
	c(x,z)=\left\{\begin{array}{ll}
		5.5, & 0<z\le 33+5\sin\frac{\pi x}{40}, \\
		7.8, & 33+5\sin\frac{\pi x}{40}<z\le 45+0.4x, \\
		7.488, & 45+0.4x<z\le 60+0.4x, \\
		8.268, & 60+0.4x<z\le 85+0.4x, \\
		7.8, & \textnormal{others}.
	\end{array}\right.
\end{equation*}	
with unit `km/s'. There are 12 randomly distributed receivers $\boeta_r=(x_r,z_r)$ on the surface $z_r=0\,km$. In Table \ref{tab:subduction_recepos}, we output their horizontal positions. This velocity models is illustrated in Figure \ref{fig:exam42_vel}. The dominant frequency of the earthquake is $f_0=2\,Hz$ and the simulating time interval $I=[0,\,55\,s]$.

\begin{table*}
    	\caption{The subduction plate model: the horizontal positions of receivers, with unit `km'.} \label{tab:subduction_recepos}
	\begin{center}\begin{tabular}{ccccccccccccc} \hline
		$r$ & $1$ & $2$ & $3$ & $4$ & $5$ & $6$ & $7$ & $8$ & $9$ & $10$ & $11$ & $12$ \\ \hline
		$x_r$ & $21$ & $33$ & $39$ & $58$ & $68$ & $74$ & $86$ & $98$ & $126$ & $132$ & $158$ & $197$ \\ \hline
    	\end{tabular}\end{center}
\end{table*}

\begin{figure*} 
	\begin{tabular}{c}
		\includegraphics[width=0.95\textwidth, height=0.12\textheight]{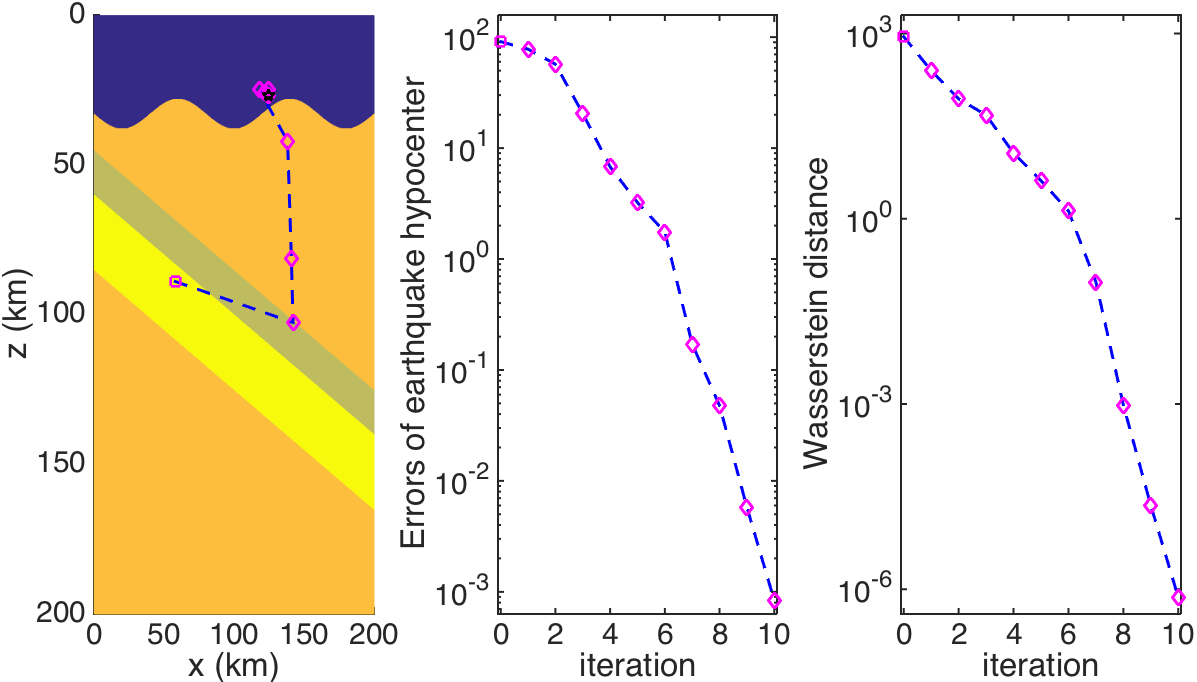} \\
		\includegraphics[width=0.95\textwidth, height=0.12\textheight]{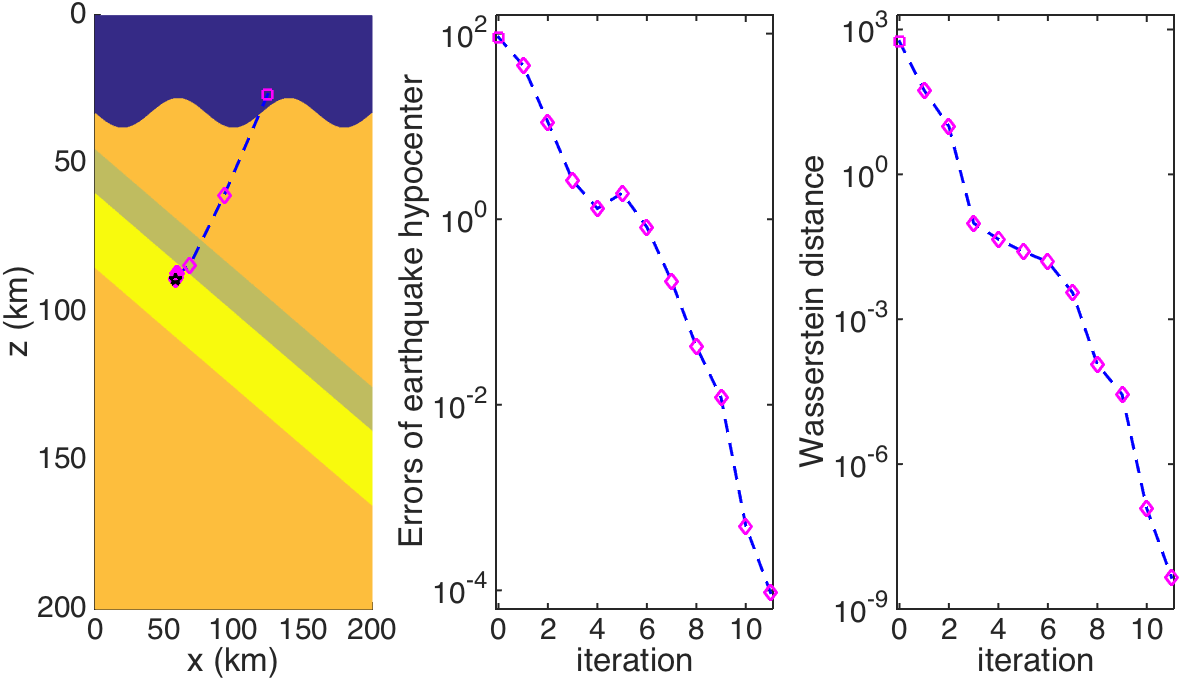}
	\end{tabular}
	\caption{Convergence history of the subduction plate model.  Up for case (i), and Down for case (ii). Left: the convergent trajectories; Mid: the absolute errors between the real and computed earthquake hypocenter with respect to iteration steps; Right: the Wasserstein distance between the real and synthetic earthquake signals with respect to iteration steps. The magenta square is the initial hypocenter, the magenta diamond denotes the hypocenter in the iterative process, and the black pentagram is the real hypocenter.} \label{fig:exam42_traj}
\end{figure*}

First, consider the ideal situation that there is no data noise. We investigate the cast that the earthquake occurs in the crust but the initial hypocenter of the earthquake is chosen in the subduction zone, and its contrary case. Their parameters are selected as follows
\begin{align*}
	(i) \; \bxi_T=(124.694\,km,\;26.762\,km),\;\tau_T=5.00\,s, \quad \bxi=(58.056\,km,\;88.985\,km),\;\tau=6.79\,s; \\
	(ii) \; \bxi_T=(58.056\,km,\;88.985\,km),\;\tau_T=6.79\,s, \quad \bxi=(124.694\,km,\;26.762\,km),\;\tau=5.00\,s.
\end{align*}
The convergent trajectories, absolute errors of the earthquake hypocenter and the Wasserstein distance are output in Figure \ref{fig:exam42_traj}. From which, we can observe nice convergence property of the new method.

\begin{figure*} 
	\begin{tabular}{cccc}
		\includegraphics[width=0.23\textwidth, height=0.12\textheight]{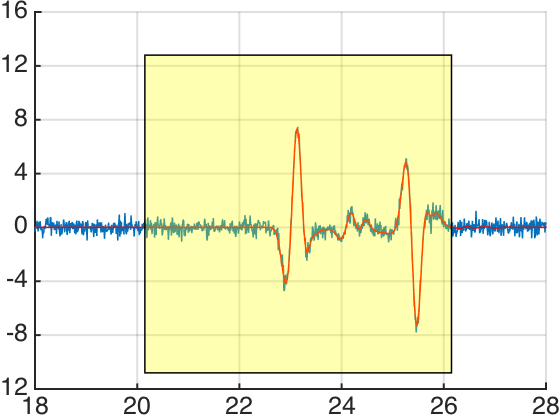} &
		\includegraphics[width=0.23\textwidth, height=0.12\textheight]{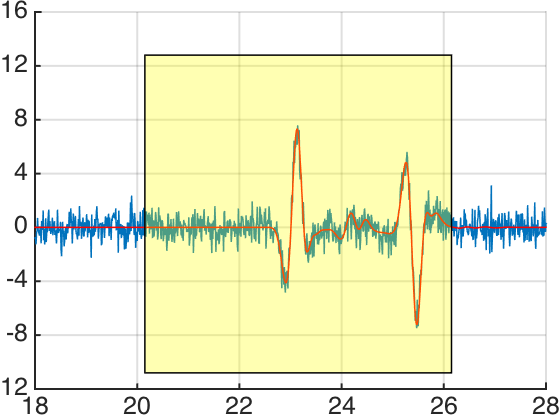} &
		\includegraphics[width=0.23\textwidth, height=0.12\textheight]{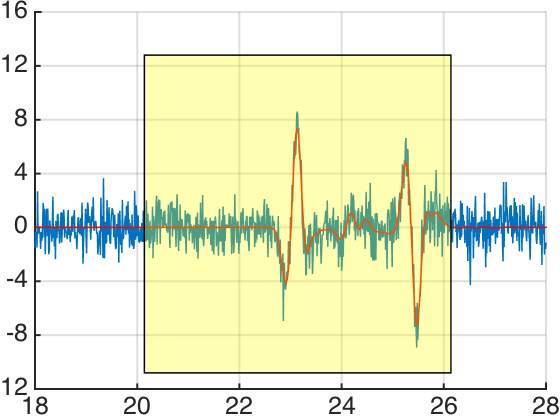} &
		\includegraphics[width=0.23\textwidth, height=0.12\textheight]{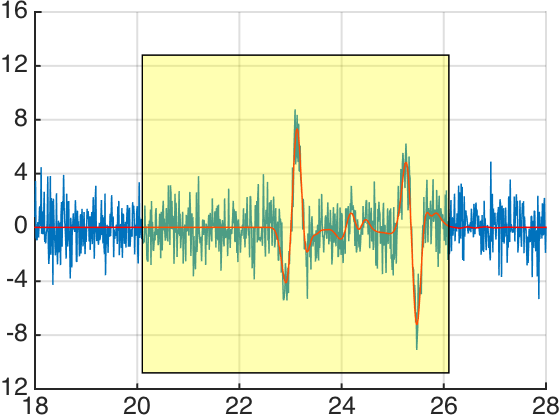} \\
		\includegraphics[width=0.23\textwidth, height=0.12\textheight]{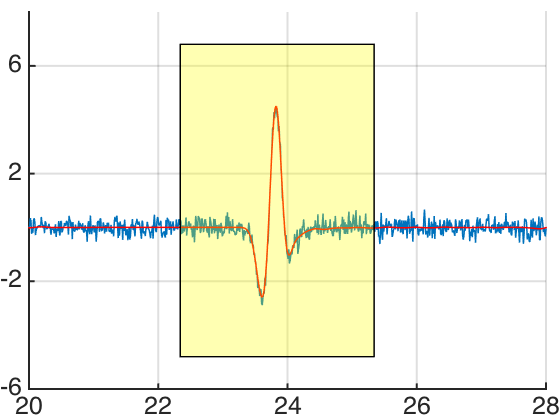} &
		\includegraphics[width=0.23\textwidth, height=0.12\textheight]{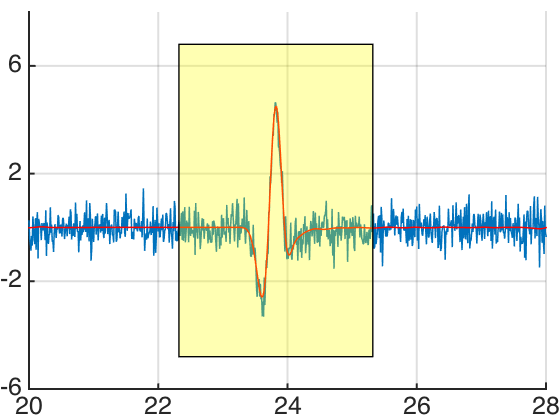} &
		\includegraphics[width=0.23\textwidth, height=0.12\textheight]{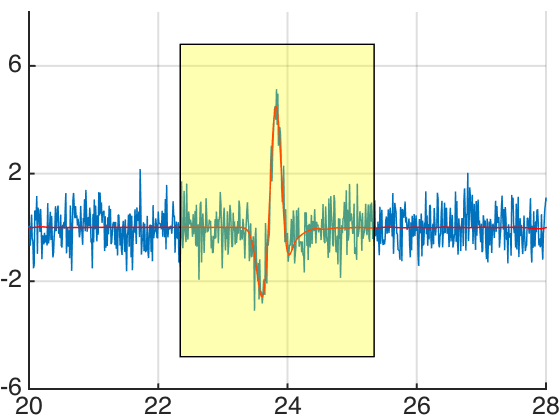} &
		\includegraphics[width=0.23\textwidth, height=0.12\textheight]{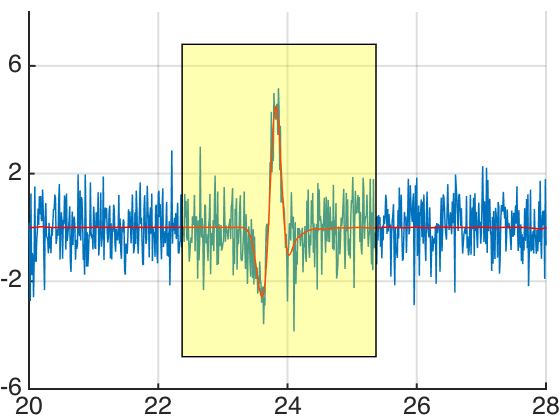}
	\end{tabular}
	\caption{Illustration of signal with noise in the subduction plate model. The signal with noise $d_r(t)$ (blue line) and the noise free signal $u(\boeta_r,t;\bxi_T,\tau_T)$ for receiver $r=4$. The horizontal axis is the time $t$. Up: parameters (i); Down: parameters (ii); From left to right, the ratio $R=5\%,\;10\%,\;15\%,\;20\%$ respectively.} \label{fig:exam42_noise_signal}
\end{figure*}

\begin{figure*} 
	\begin{tabular}{c}
		\includegraphics[width=0.95\textwidth, height=0.10\textheight]{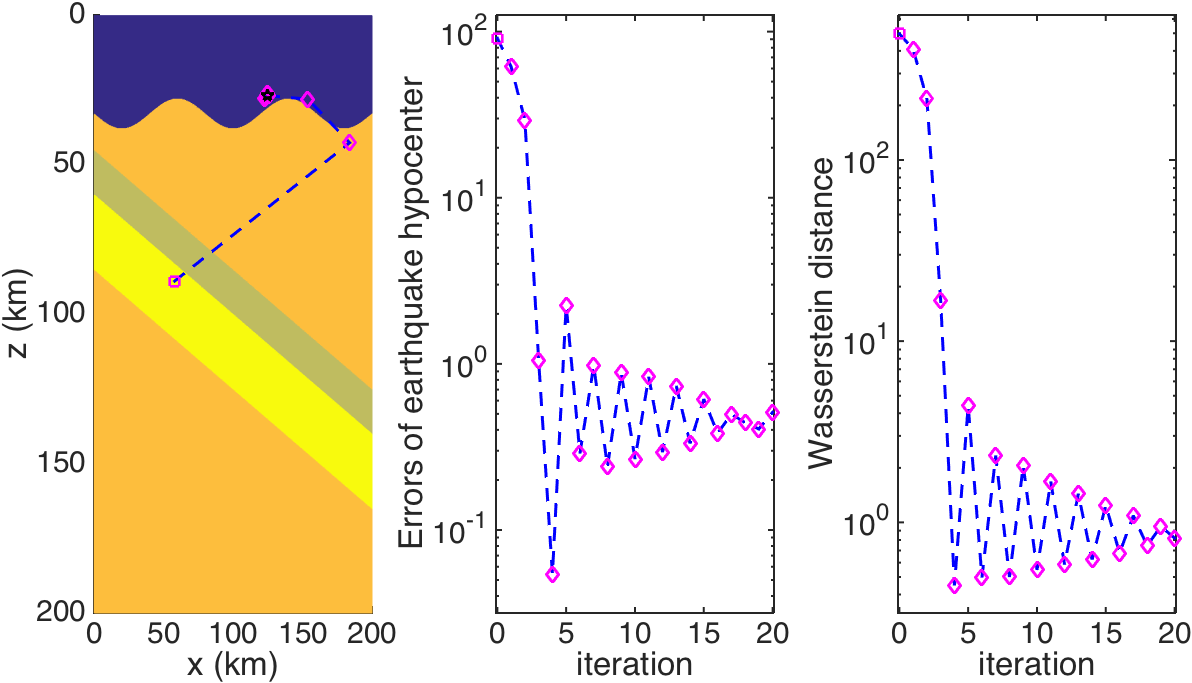} \\
		\includegraphics[width=0.95\textwidth, height=0.10\textheight]{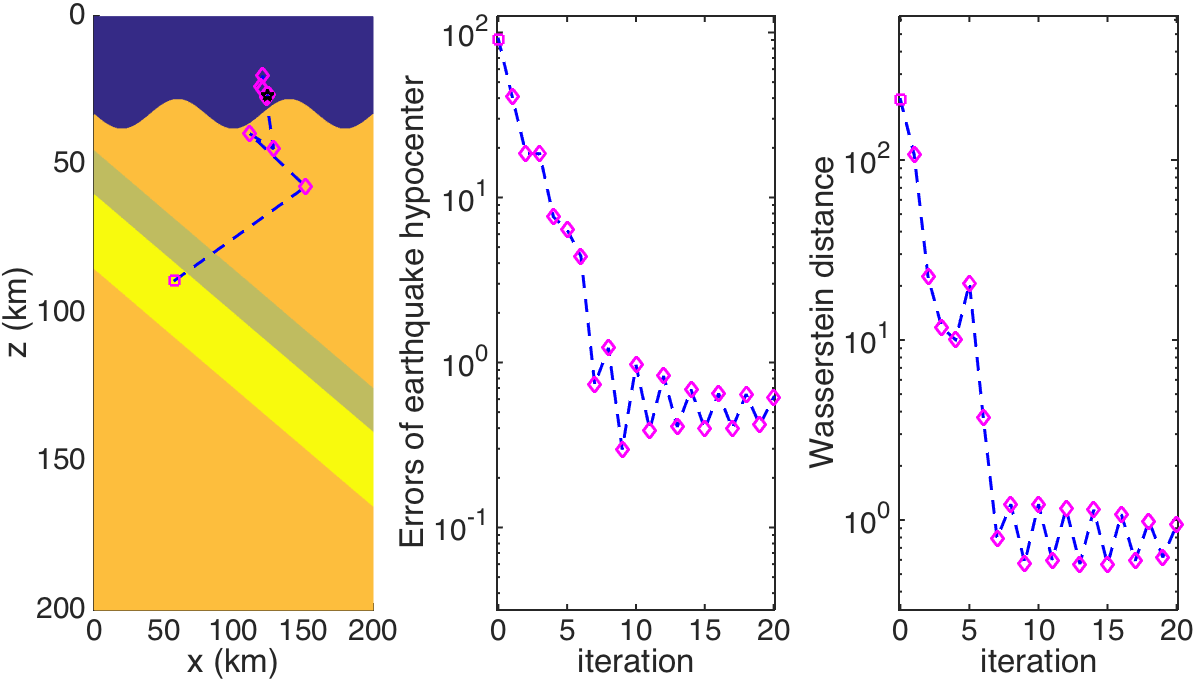} \\
		\includegraphics[width=0.95\textwidth, height=0.10\textheight]{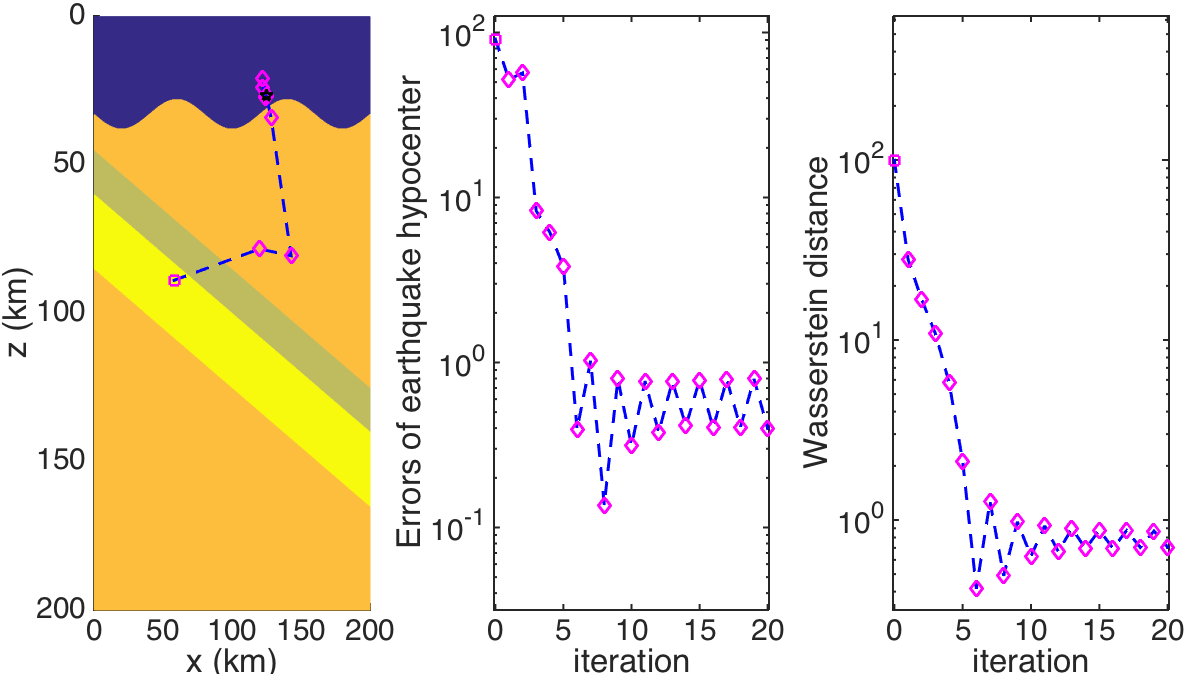} \\
		\includegraphics[width=0.95\textwidth, height=0.10\textheight]{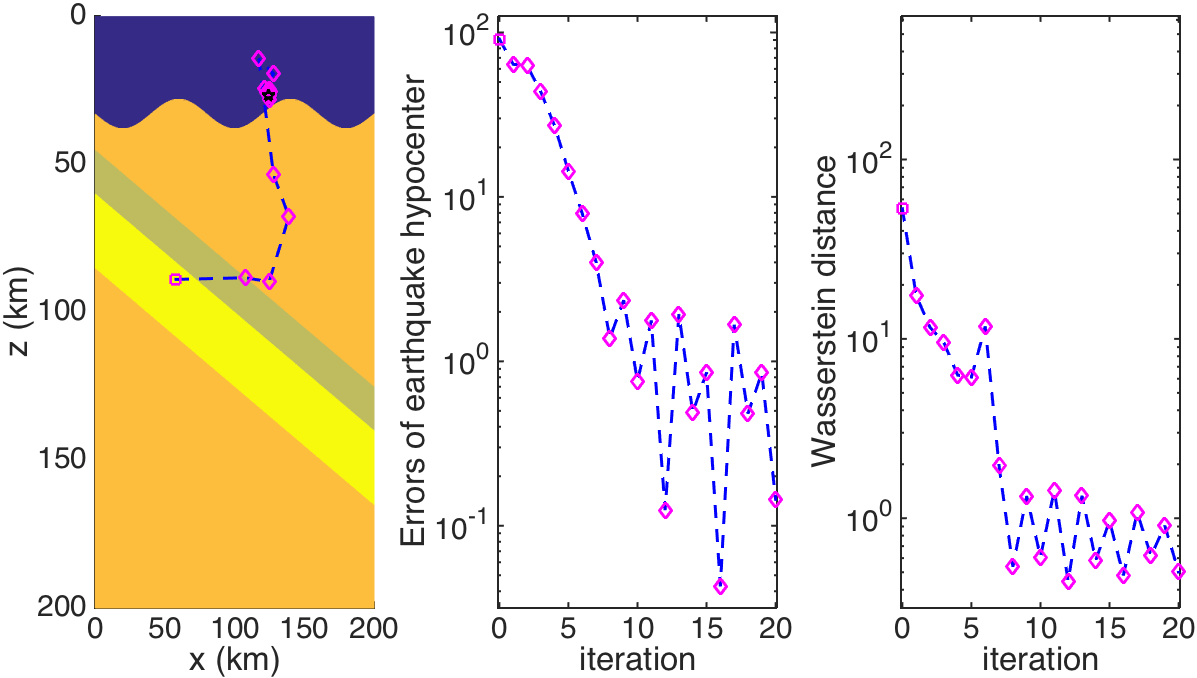}
	\end{tabular}
	\caption{Convergence history of the subduction plate model with noise data, case (i). From up to bottom, the the ratio $R=5\%,\;10\%,\;15\%,\;20\%$ respectively. Left: the convergent trajectories; Mid: the absolute errors between the real and computed earthquake hypocenter with respect to iteration steps; Right: the Wasserstein distance between the real and synthetic earthquake signals with respect to iteration steps. The magenta square is the initial hypocenter, the magenta diamond denotes the hypocenter in the iterative process, and the black pentagram is the real hypocenter.} \label{fig:exam42_noise1_traj}
\end{figure*}

\begin{figure*} 
	\begin{tabular}{c}
		\includegraphics[width=0.95\textwidth, height=0.10\textheight]{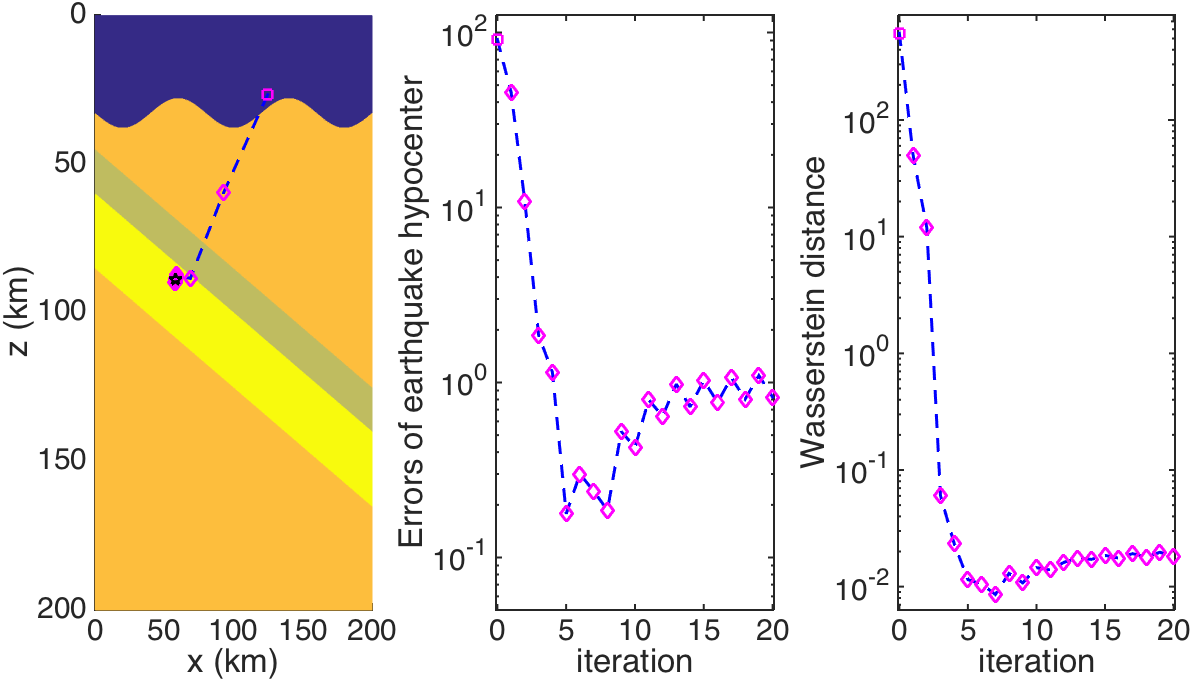} \\
		\includegraphics[width=0.95\textwidth, height=0.10\textheight]{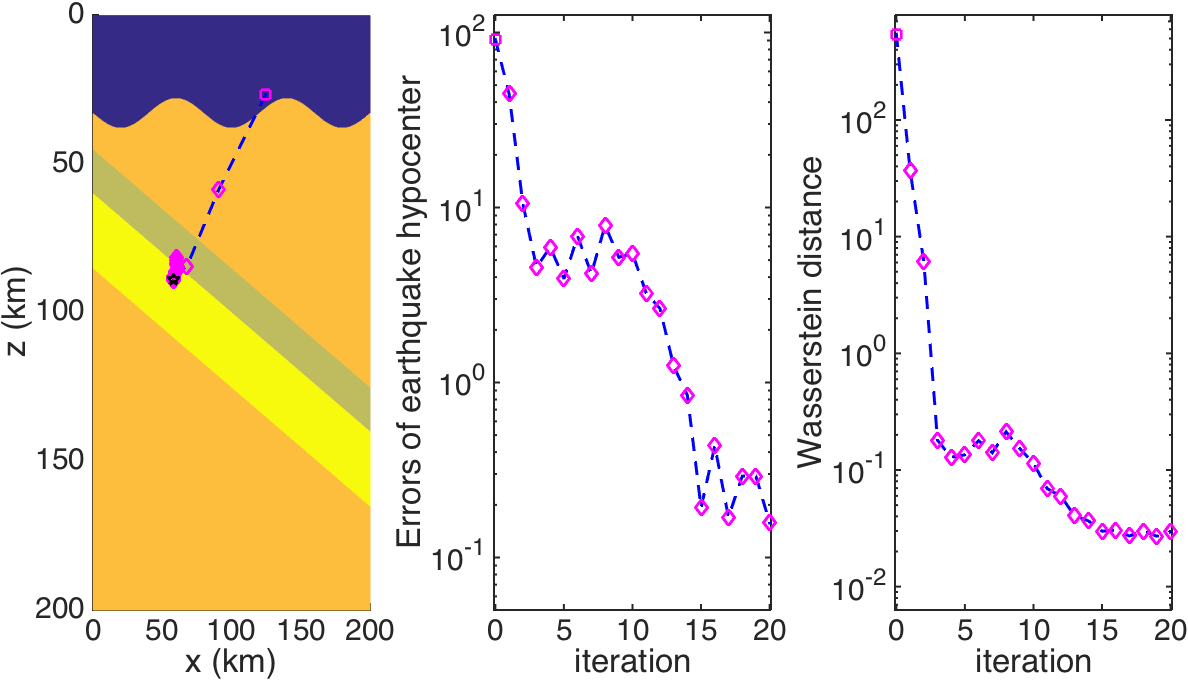} \\
		\includegraphics[width=0.95\textwidth, height=0.10\textheight]{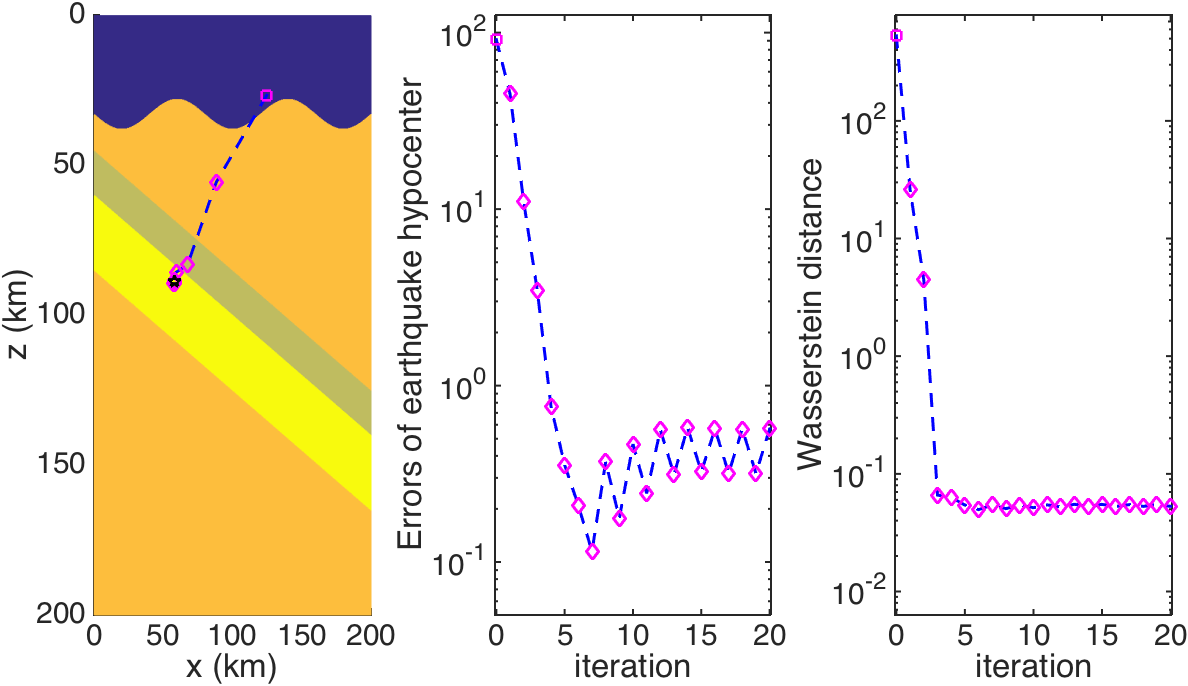} \\
		\includegraphics[width=0.95\textwidth, height=0.10\textheight]{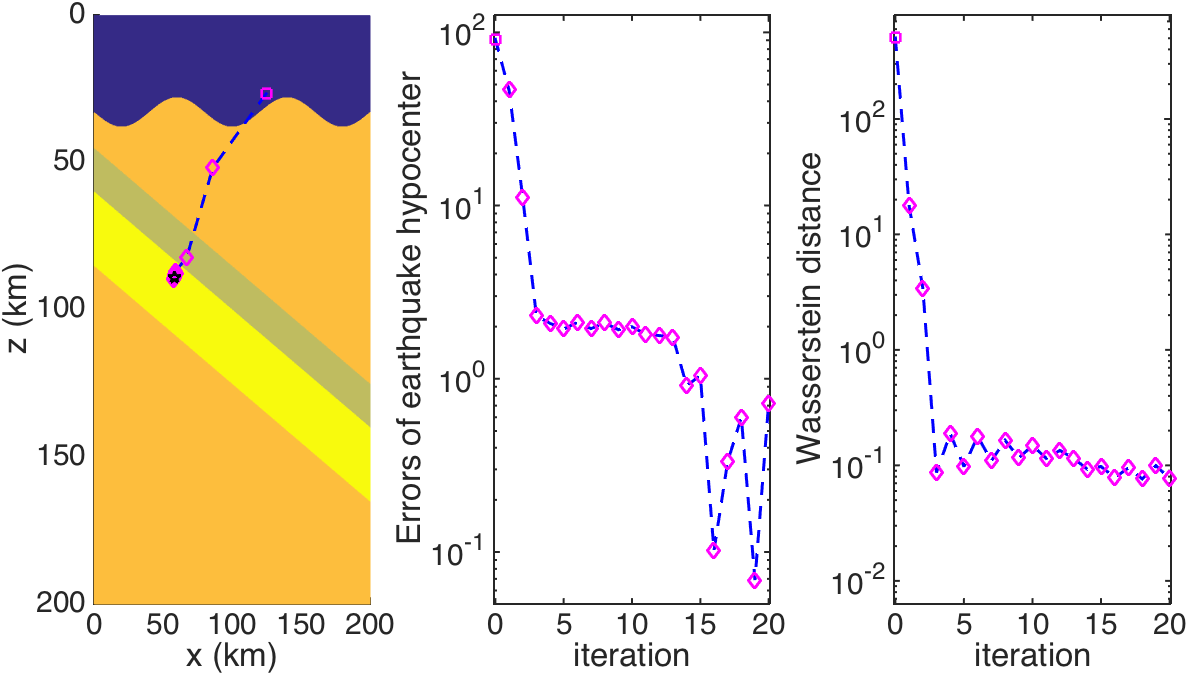}
	\end{tabular}
	\caption{Convergence history of the subduction plate model with noise data, case (ii). From up to bottom, the the ratio $R=5\%,\;10\%,\;15\%,\;20\%$ respectively. Left: the convergent trajectories; Mid: the absolute errors between the real and computed earthquake hypocenter with respect to iteration steps; Right: the Wasserstein distance between the real and synthetic earthquake signals with respect to iteration steps. The magenta square is the initial hypocenter, the magenta diamond denotes the hypocenter in the iterative process, and the black pentagram is the real hypocenter.} \label{fig:exam42_noise2_traj}
\end{figure*}

We next consider the signal containing noise. We select the same parameters (i) and (ii). The noise is added to the real earthquake signals in the same way as in Subsection \ref{subsec:twolayer}. In Figure \ref{fig:exam42_noise_signal}, the real earthquake signal with noise $d_r(t)$ and the noise free signal $u(\boeta_r,t;\bxi_T,\tau_T)$ are presented. In order to reduce the impact of noise, it is necessary and reasonable to select a time window that contains the main part of $u(\boeta_r,t;\bxi_T,\tau_T)$. Moreover, the technique proposed in Remark \ref{rem:noise} and the modified LMF method (Algorithm \ref{alg:LMF_nd}) are also applied here. In Figure \ref{fig:exam42_noise1_traj}-\ref{fig:exam42_noise2_traj}, we output the convergent history. As discussed in the previous Subsection, the location errors and the misfit functions oscillate due to the noise effect. Thus, it is impossible to get accurate location results. A reasonable choice is to select the iteration step $k_*$ which corresponds to the smallest value of the misfit function. These values are output in Table \ref{tab:exam42_noise1}-\ref{tab:exam42_noise2}. From which, we can see the location results are more satisfying by the new model and method.

\begin{table*}
    	\caption{The subduction plate mode with noise data, case (i). The smallest misfit value, the corresponding iteration step $k_*$ and the location error.} \label{tab:exam42_noise1}
	\begin{center}\begin{tabular}{c|ccc} \hline
		 $R$ & $k_*$ & The misfit value & The location error (km) \\ \hline
		 $5\%$ & $4$ & $4.49\times 10^{-1}$ & $5.41\times 10^{-2}$ \\ 
		 $10\%$ & $15$ & $5.67\times 10^{-1}$ & $3.96\times 10^{-1}$ \\ 
		 $15\%$ & $6$ & $4.15\times 10^{-1}$ & $3.93\times 10^{-1}$ \\ 
		 $20\%$ & $12$ & $4.45\times 10^{-1}$ & $1.24\times 10^{-1}$ \\ \hline
    	\end{tabular}\end{center}
\end{table*}

\begin{table*}
    	\caption{The subduction plate mode with noise data, case (ii). The smallest misfit value, the corresponding iteration step $k_*$ and the location error.} \label{tab:exam42_noise2}
	\begin{center}\begin{tabular}{c|ccc} \hline
		 $R$ & $k_*$ & The misfit value & The location error (km) \\ \hline
		 $5\%$ & $7$ & $8.56\times 10^{-3}$ & $2.40\times 10^{-1}$ \\ 
		 $10\%$ & $19$ & $2.70\times 10^{-2}$ & $2.92\times 10^{-1}$ \\ 
		 $15\%$ & $6$ & $4.96\times 10^{-2}$ & $2.10\times 10^{-1}$ \\ 
		 $20\%$ & $18$ & $7.64\times 10^{-2}$ & $5.99\times 10^{-1}$ \\ \hline
    	\end{tabular}\end{center}
\end{table*}

\section{Conclusion} \label{sec:con}
In this paper, we apply the quadratic Wasserstein metric to the earthquake location problem. The numerical evidence suggests that the convexity of the misfit function with respect to the earthquake hypocenter and origin time, which is based on the quadratic Wasserstein metric, is much better than the one based on the $L^2$ metric. This makes it possible to accurately locate the earthquakes even starting from very far initial values. Besides, since the misfit function is close to the quadratic function, the LMF method could be a good choice to solve the resulted optimization problem. According to our numerical tests, the LMF method has obvious advantages over the GN method and the BFGS method. When the original signal is affected by noise, we make a little modification of the quadratic Wasserstein metric based misfit function and the LMF method. In all the numerical experiments, the location errors are smaller than $1\,$km. These location results, according to our best knowledge, are pretty good.

We also need to point out that, since both the real and synthetic signals have been normalized, it is not possible to determine the amplitude of the seismogram at source. To deal with this difficulty, we may need to introduce the unbalanced optimal transport theory \cite{ChPeScVi:16, ChPeScVi:17}. Moreover, the techniques developed in this paper may be applicable to the inversion of many micro-earthquakes \cite{LeSt:81, PrGe:88}. We are currently working on these interesting topics and hope to report these in subsequent papers.

\section*{Acknowledgement}
H. Wu was supported by the NSFC project 11101236, NSFC project 91330203 and SRF for ROCS, SEM. D.H. Yang was supported by NSFC project 41230210 and NSFC project 41390452. The authors are grateful to Prof. Shi Jin for his inspiration and helpful suggestions and discussions that greatly improve the presentation. Hao Wu would like to acknowledge Prof. Bj\"orn Engquist and Prof. Brittany D. Froese for their valuable comments.


\begin{thebibliography}{}
	\bibitem{AmGi:13}
		L. Ambrosio and N. Gigli, A user guide to optimal transport, in \textit{Modelling and Optimisation of Flows on Networks}, pp. 1-155, Springer, 2013.
		
	\bibitem{ArChBo:17}
		M. Arjovsky, S. Chintala and L. Bottou, Wasserstein GAN, \textit{arXiv:1701.07875}, 2017.

	\bibitem{AkRi:80}
		K. Aki and P.G. Richards, \textit{Quantitative Seismology: Theory and Methods volume II}, W.H. Freeman \& Co (Sd), 1980.
		
	\bibitem{ChPeScVi:16}
		L. Chizat, G. Peyr\'e, B. Schmitzer and F.X. Vialard, An Interpolating Distance Between Optimal Transport and Fisher-Rao Metrics, \textit{Found. Comput. Math.}, 2016. https://doi.org/10.1007/s10208-016-9331-y
		
	\bibitem{ChPeScVi:17}
		L. Chizat, G. Peyr\'e, B. Schmitzer and F.X. Vialard, Scaling Algorithms for Unbalanced Optimal Transport Problems, \textit{arXiv:1607.05816v2}, 2017.	
				
	\bibitem{Da:86}
		M.A. Dablain, The application of high-order differencing to the scalar wave equation, \textit{Geophysics}, \textbf{51}(1), 54-66, 1986.			

	\bibitem{EnFr:14}
    		B. Engquist and B.D. Froese, Application of the Wasserstein metric to seismic signals, \textit{Commun. Math. Sci.}, \textbf{12}(5), 979-988, 2014.

	\bibitem{EnFrYa:16}
    		B. Engquist, B.D. Froese and Y.N. Yang, Optimal transport for seismic full waveform inversion, \textit{Commun. Math. Sci.}, \textbf{14}(8), 2309-2330, 2016.

	\bibitem{EnRu:03}
    		B. Engquist and O. Runborg, Computational high frequency wave propagation, \textit{Acta Numer.}, \textbf{12}, 181-266, 2003.	
		
	\bibitem{Fl:91}
		R. Fletcher, \textit{Practical Methods of  Optimization}, Second edition, Chichester: John Wiley and Sons, 1991.

	\bibitem{Ge:03}
    		M.C. Ge, Analysis of source location algorithms Part I: Overview and non-iterative methods, \textit{J. Acoust. Emiss.}, \textbf{21}, 14-28. 2003.

	\bibitem{Ge:03b}
    		M.C. Ge, Analysis of source location algorithms Part II: Iterative methods, \textit{J. Acoust. Emiss.}, \textbf{21}, 29-51, 2003.

	\bibitem{Ge:12}
    		L. Geiger, Probability method for the determination of earthquake epicenters from the arrival time only, \textit{Bull. St. Louis Univ.}, \textbf{8}, 60-71, 1912.

	\bibitem{HyKaOj:01}
		A. Hyv\"arinen, J. Karhunen and E. Oja, \textit{Independent Component Analysis}, John Wiley \& Sons, Inc., 2001.

	\bibitem{JiWuYa:08}
    		S. Jin, H. Wu and X. Yang, Gaussian Beam Methods for the Schr\"odinger Equation in the Semi-classical Regime: Lagrangian and Eulerian Formulations, \textit{Commun. Math. Sci.}, \textbf{6}(4), 995-1020, 2008.
		
	\bibitem{KaRu:58}
		L.V. Kantorovich and G.S. Rubinshtein, On a space of totally additive functions, \textit{Vestn. Leningrad. Univ.}, \textbf{13}(7), 52-59, 1958.
				
	\bibitem{KiLiTr:11}
		Y.H. Kim, Q.Y. Liu and J. Tromp, Adjoint centroid-moment tensor inversions, \textit{Geophys. J. Int.}, \textbf{186}, 264-278, 2011.
		
	\bibitem{KoTr:03}
		D. Komatitsch and J. Tromp, A perfectly matched layer absorbing boundary condition for the second-order seismic wave equation, \textit{Geophys. J. Int.}, \textbf{154}, 146-153, 2003.
		
	\bibitem{LeSt:81}
		W.H.K. Lee and S.W. Stewart, \textit{Principles and Applications of Microearthquake Networks}, Academic Press, 1981.

	\bibitem{LeSt:81}
		W.H.K. Lee and S.W. Stewart, \textit{Principles and Applications of Microearthquake Networks}, Academic Press, 1981.

	\bibitem{Le:44}
		K. Levenberg, A Method for the Solution of Certain Problems in Least Square, \textit{Quart. Appl. Math.}, \textbf{2}, 164-168, 1944.
		
	\bibitem{LiYaWuMa:17}
		J.S. Li, D.H. Yang, H. Wu and X. Ma, A low-dispersive method using the high-order stereo-modelling operator for solving 2-D wave equations, \textit{Geophys. J. Int.}, \textbf{210}, 1938-1964, 2017.
		
	\bibitem{LiPoKoTr:04}
		Q.Y. Liu, J. Polet, D. Komatitsch and J. Tromp, Spectral-Element Moment Tensor Inversion for Earthquakes in Southern California, \textit{Bull. seism. Soc. Am.}, \textbf{94}(5), 1748-1761, 2004.
		
	\bibitem{Ma:15}
		R. Madariaga, Seismic Source Theory, in \textit{Treatise on Geophysics (Second Edition)}, pp. 51-71, ed. Gerald, S., Elsevier B.V., 2015.
		
	\bibitem{MaNiTi:04}
		K. Madsen, H.B. Nielsen and O. Tingleff, \textit{Methods for Non-Linear Least Squares Problems (2nd ed.)}, Informatics and Mathematical Modelling, Technical University of Denmark, 2004.

	\bibitem{Ma:63}
		D. Marquardt, An Algorithm for Least Squares Estimation on Nonlinear Parameters, \textit{SIAM J. Appl. Math.}, \textbf{11}, 431-441, 1963.
		
	\bibitem{MeBrMeOuVi:16}
		L. M\'etivier, R. Brossier, Q. M\'erigot, E. Oudet and J. Virieux, Measuring the misfit between seismograms using an optimal transport distance: application to full waveform inversion, \textit{Geophys. J. Int.}, \textbf{205}, 345-377, 2016.
		
	\bibitem{MeBrMeOuVi:16b}
		L. M\'etivier, R. Brossier, Q. M\'erigot, E. Oudet and J. Virieux, An optimal transport approach for seismic tomography: application to 3D full waveform inversion, \textit{Inverse Probl.}, \textbf{32}, 115008, 2016.

	\bibitem{Mo:81}
		G. Monge, M\'emoire sur la th\'eorie des d\'eblais et de remblais, \textit{Histoire de l'acad\'emie royale des sciences de paris, avec les M\'emoires de Math\'ematique et de Physique pour la mme ann\'ee}, pp. 666–704, 1781.

	\bibitem{PrGe:88}
		A.F. Prugger and D.J. Gendzwill, Microearthquake location: A nonlinear approach that makes use of a simplex stepping procedure, \textit{Bull. seism. Soc. Am.}, \textbf{78}, 799-815, 1988.
		
	\bibitem{RaPoFi:10}
    		N. Rawlinson, S. Pozgay and S. Fishwick, Seismic tomography: A window into deep Earth, \textit{Phys. Earth Planet. Inter.}, \textbf{178}, 101-135, 2010.

	\bibitem{RuToGu:98}
		Y. Rubner, C. Tomasi and L. J. Guibas, A metric for distributions with applications to image databases, in \textit{IEEE International Conference on Computer Vision}, pp. 59-66, 1998.
				
	\bibitem{Sa:15}
		F. Santambrogio, \textit{Optimal Transport for Applied Mathematicians: Calculus of Variations, PDEs and Modeling}, Progress in Nonlinear Differential Equations and Their Applications, Birkh\"auser, 2015.
				
	\bibitem{SaLoZo:08}
		C. Satriano, A. Lomax and A. Zollo, Real-Time Evolutionary Earthquake Location for Seismic Early Warning, \textit{Bull. seism. Soc. Am.}, \textbf{98}(3), 1482-1494, 2008.
				
	\bibitem{Th:85}
		C.H. Thurber, Nonlinear earthquake location: Theory and examples, \textit{Bull. seism. Soc. Am.}, \textbf{75}(3), 779-790, 1985.

	\bibitem{Th:14}
		C.H. Thurber, Earthquake, location techniques, in \textit{Encyclopedia of Earth Sciences Series}, pp. 201-207, ed. Gupta, H.K., Springer, 2014.
	
	\bibitem{ToYaLiYaHa:16}
		P. Tong, D.H. Yang, Q.Y. Liu, X. Yang and J. Harris, Acoustic wave-equation-based earthquake location, \textit{Geophys. J. Int.}, \textbf{205}(1), 464-478, 2016.

	\bibitem{ToZhYa:11}
		P. Tong, D.P. Zhao and D.H. Yang, Tomography of the 1995 Kobe earthquake area: comparison of finite-frequency and ray approaches, \textit{Geophys. J. Int.}, \textbf{187}, 278-302, 2011.					
	\bibitem{Vi:03}
		C. Villani, \textit{Topics in Optimal Transportation}, Graduate Studies in Mathematics, American Mathematical Society, 2003.

	\bibitem{Vi:08}
		C. Villani, \textit{Optimal Transport: Old and New}, Springer Science \& Business Media, 2008.

	\bibitem{WaEl:00}
		F. Waldhauser and W.L. Ellsworth, A double-difference earthquake location algorithm: Method and application to the northern Hayward Fault, California, \textit{Bull. seism. Soc. Am.}, \textbf{90}(6), 1353-1368, 2000.
		
	\bibitem{We:08}
    		X. Wen, High Order Numerical Quadratures to One Dimensional Delta Function Integrals, \textit{SIAM J. Sci. Comput.}, \textbf{30}(4), 1825-1846, 2008.		

	\bibitem{WuChHuYa:16}
    		H. Wu, J. Chen, X.Y. Huang and D.H. Yang, A new earthquake location method based on the waveform inversion, \textit{Commun. Comput. Phys.}, \textbf{23}(1), 118-141, 2018.
					
	\bibitem{WuChJiJiYa:17}
		H. Wu, J. Chen, H. Jing, P. Tong and D.H. Yang, The auxiliary function method for waveform based earthquake location, \textit{arXiv:1706.05551}, 2017.
				
	\bibitem{WuYa:13}
    		H. Wu and X. Yang, Eulerian Gaussian beam method for high frequency wave propagation in the reduced momentum space, \textit{Wave Motion}, \textbf{50}(6), 1036-1049, 2013.

	\bibitem{YaLuWuPe:04}
		D.H. Yang, M. Lu, R.S. Wu and J.M. Peng, An Optimal Nearly Analytic Discrete Method for 2D Acoustic and Elastic Wave Equations, \textit{Bull. seism. Soc. Am.}, \textbf{94}(5), 1982-1991, 2004.

	\bibitem{YaEnSuFr:16}
    		Y.N. Yang, B. Engquist, J.Z. Sun and B.D. Froese, Application of Optimal transport and the quadratic Wasserstein metric to Full-Waveform-Inversion, \textit{arXiv:1612.05075}, 2016.			
\end{thebibliography}
\end{document}